\documentclass[11pt]{article}

\usepackage{graphicx}
\usepackage{marvosym}
\usepackage{amsmath}
\usepackage{amssymb}
\usepackage{amsthm}
\usepackage{tikz}
\usepackage{tikz-cd}
\usepackage[left = 1.3in, right = 1.3in]{geometry}
\usepackage[utf8]{inputenc}
\usepackage[english]{babel}
\usepackage[maxbibnames=99]{biblatex}
\usepackage{hyperref}
\usepackage{todonotes}
\usepackage{xcolor}

\newcommand{\complex}{\mathbb{C}}

\newcommand{\reals}{\mathbb{R}}
\newcommand{\intgr}{\mathbb{Z}}
\newcommand{\ntrl}{\mathbb{N}}
\newcommand{\rtnl}{\mathbb{Q}}

\renewcommand{\tilde}{\widetilde}
\renewcommand{\subset}{\subseteq}
\renewcommand{\phi}{\varphi}
\renewcommand{\epsilon}{\varepsilon}

\renewcommand{\hat}{\widehat}

\newcommand{\tr}{\text{tr}}

\renewcommand{\exp}{\text{exp}}

\renewcommand{\H}{\mathbb{H}}
\renewcommand{\exp}{\textrm{exp}}

\newcommand{\M}{\mathcal{M}}

\newcommand{\Res}{\mathrm{Res}}
\newcommand{\ep}{\epsilon}
\newcommand{\re}{\mathrm{Re}}
\renewcommand{\L}{\mathcal{L}}

\newtheorem{theorem}{Theorem}[section]
\newtheorem{lemma}[theorem]{Lemma}
\newtheorem{corollary}[theorem]{Corollary}

\theoremstyle{definition}
\newtheorem{definition}[theorem]{Definition}
\theoremstyle{remark}
\newtheorem{remark}[theorem]{Remark}
\newtheorem{example}{Example}[section]

\title{Dynamical Zeta Functions and Resonance Chains for Infinite-Area Hyperbolic Surfaces with Large Funnel Widths}

\author{Henry Talbott}

\begin{document}

\maketitle

\begin{abstract}

We quantitatively relate the resonance sets of topologically finite infinite-area hyperbolic surfaces with no cusps to the resonance sets of certain metric graphs via the spine graph construction. In particular, we prove the existence of approximate resonance chains in resonance sets of these surfaces in the long-boundary-length regime. Our results are similar in spirit to those obtained in recent independent work by Li-Matheus-Pan-Tao, although our perspective and hypotheses are somewhat different. Our results also generalize older results obtained for three-funneled spheres by Weich. We primarily make use of transfer operators for holomorphic iterated function schemes, along with certain geometric bounds.

\end{abstract}

\tableofcontents

\section{Introduction}

Recall that $\M_{g,n}(L_1,\dots,L_n)$ is the moduli space of hyperbolic surfaces with genus $g$, $n$ boundary components, and prescribed boundary lengths $L_1,\dots,L_n$. It is naturally also the moduli space of hyperbolic surfaces with genus $g$, $n$ infinite-area funnels, and prescribed funnel widths $L_1,\dots,L_n$. Here, the `width' of a funnel is the length of the unique geodesic separating that funnel from the rest of the surface. To realize this equivalence, take each boundary geodesic of a surface with boundary, and glue on the unique funnel with the correct width. We will think of $\M_{g,n}(L_1,\dots,L_n)$ using this second identification. 

For any such $X\in\M_{g,n}(L_1,\dots,L_n)$, we may define its Laplacian operator $\Delta_X$. The \emph{resolvent} of $X$ is the following family of operators, varying according to a parameter $s\in\complex$:
\begin{equation}R_X(s)=(\Delta_X-s(s-1))^{-1}.\end{equation}

Formally, the above definition of $R_X(s)$ is only well-defined for $\re(s)\geq 1/2$. However, we have:

\begin{theorem}[Mazzeo-Melrose \cite{MM87}, Guillop{\'e}-Zworski, Theorem 1 \cite{GZ95}]
\label{resolventconst}

The resolvent $R_X(s)$ defined above extends to a meromorphic family of operators
$$R_X(s):L^2_{\textrm{comp}}(X)\to H^2_{\textrm{loc}}(X)$$
with poles of finite rank. 
	
\end{theorem}

We may therefore define the \emph{resonance set} of $X$:

$$\Res(X)=\{s\in\complex:s\textrm{ is a pole of }R_X(s)\}$$
where poles are counted with multiplicity equal to their rank. 

These resonances can be thought of as generalized Laplacian eigenvalues, and play an important role in the study of infinite-area hyperbolic surfaces. For a reference work on this theory, see \cite{DB16}.

Many significant results have been proven about the general distributional properties of these resonances; a survey is provided in \cite{SN11}. Furthermore, questions about resonances link hyperbolic geometry to several other fields, including number theory  \cite{BGS11} and the study of chaotic quantum systems in theoretical physics \cite{BKPSWZ12,BKPSWZ13}.

In this paper, we will show that the resonances of an infinite-area surface $X$ can be related to the zeroes of a certain dynamical zeta function derived from a different geometric object,  the \emph{spine graph} of $X$. While this relation is only approximate, for {`typical'} surfaces, the {error} terms decay exponentially in the minimum length of the funnel widths of $X$. 

Specifically, for any $\M_{g,n}(L_1,\dots,L_n)$, there exists a `surface-to-spine' homeomorphism
\begin{equation}\Phi:\M_{g,n}(L_1,\dots,L_n)\to MRG_{g,n}(L_1,\dots,L_n)\end{equation}
where $MRG_{g,n}(L_1,\dots,L_n)$ is the corresponding moduli space of \emph{metric ribbon graphs} or \emph{fatgraphs}. Metric ribbon graphs are metric graphs with extra vertex orientation information (see \cite{ABCGLW21}). {Note that here both spaces carry the standard topology; the topology of $\M_{g,n}(L_1,...,L_n)$ may be defined through Fenchel-Nielsen coordinate charts, whereas the topology on $MRG_{g,n}(L_1,...,L_n)$ may be defined via a natural cellularization, as in \cite{ABCGLW21}.} This homeomorphism was first defined by Bowditch-Epstein \cite{BE88}, and has been used in the study of hyperbolic surfaces by Do, Mondello, and others \cite{ND10,GM09,HT25}. 

To set up the resonance correspondence, we first observe that associated to any $X\in\M_{g,n}(L_1,\dots,L_n)$ is a \emph{dynamical zeta function} $d_X(s,z)$, an entire function in both complex variables satisfying
$$\Res(X)=\{s\in\complex:d_X(s,1)=0\}\setminus(-\intgr_{\geq 1}).$$

For any $\Gamma\in MRG_{g,n}(L_1,\dots,L_n)$, there likewise exists a metric graph dynamical zeta function $\tilde{d}_\Gamma(s,z)$, and we may define
$$\tilde{\Res}(\Gamma)=\{s\in\complex:\tilde{d}_\Gamma(s,1)=0\}.$$
{The definition of $\tilde{d}_\Gamma(s,z)$ can be found in Subsection \ref{zetafunctions}, while the definition of $d_X(s,z)$ is delayed until Lemma \ref{selbergequality}.} We will argue that $\tilde{\Res}(\Phi(X))$ converges to $\Res(X)$ under the correct limit and the correct notion of convergence.

\subsection{Main Results}

Graph resonances are scaling-invariant, in the sense that for $\alpha\in\reals_{>0}$, $\tilde{\Res}(\alpha\Gamma)=\frac{1}{\alpha}\tilde{\Res}(\Gamma)$. {Here, $\alpha\Gamma$ is the metric ribbon graph obtained by scaling each edge length of $\Gamma$ by $\alpha$.} It therefore makes sense to define scaled resonance sets
$$\Res(X,L)=\{s\in\complex:\frac{s}{L}\in\Res(X)\},$$
$$\tilde{\Res}(\Gamma,L)=\{s\in\complex:\frac{s}{L}{\in}\tilde{\Res}(\Gamma)\}=\tilde{\Res}(\frac{1}{L}\Gamma).$$
We also need a notion of an `$\eta$-great surface', {where $\eta$ is some positive real number}. This condition is almost that $\Phi(X)$ is trivalent and all edges in $\Phi(X)$ have length at least $2\eta$; see Definition \ref{gooddef} for a precise statement.

Lastly, we need the following definition:

\begin{definition}

For a moduli space $\M_{g,n}(L_1,\dots,L_n)$ and given parameters $L$ and $C\geq 1$, we say that the $L_h$ (with $1\leq h\leq n$) are $C$-bounded by $L$ if they satisfy $L_h\leq CL$ for all $L_h$, and that the $L_h$ are $C$-controlled by $L$ if they satisfy
$$\frac{1}{C}L\leq L_h\leq CL$$
for all $L_h$.

\end{definition}

Note that the $\eta$-greatness condition is not so strict. Take $\M_{g,n}(L_1,\dots,L_n)$ and assume the $L_h$ are $C$-controlled by $L$. Also set $\eta=L^\epsilon$ for some $0<\epsilon<1$. Then by Corollary 4.13 in \cite{HT25}, the Weil-Petersson probability that $X\in\M_{g,n}(L_1,\dots,L_n)$ is \emph{not} $\eta$-great is $O_{g,n,C}(1/L^{2-2\epsilon})$.

Our main results are as follows. We take a surface $X\in\M_{g,n}(L_1,...,L_n)$ and the corresponding metric ribbon graph $\Gamma=\Phi(X)\in MRG_{g,n}(L_1,...L_n)$, and associate respective dynamical zeta functions $d_X(s,z)$ and $\tilde{d}_\Gamma(s,z)$ to both objects. We then show that if all $L_h$ are large, these two zeta functions approximately coincide on compact sets. Specifically:

\begin{theorem}
\label{mainresult}

Fix $X\in\M_{g,n}(L_1,\dots,L_n)$, and assume that $X$ is $\eta$-great and that the $L_h$ are $C$-bounded by $L$. Set $\Gamma=\Phi(X)$, and fix $A\in\reals_{>0}$. Also let $\epsilon>0$. Assume that $\eta$ is sufficiently large in terms of $g$, $n$, $C$, $A$, {and $\epsilon$}. Then there exist dynamical zeta functions $d_X(s,z)$ and $\tilde{d}_\Gamma(s,z)$, entire in both variables, such that for all $s,z\in\complex$ satisfying $|z|\leq A$ and $|s|\leq\frac{A}{L}$, we have
\begin{equation}|d_X(s,z)-\tilde{d}_\Gamma(s,z)|=O_{g,n,C,A,{\epsilon}}(e^{-{(1-\epsilon)\eta}})\end{equation}
where the implicit constant depends only on $g$, $n$, $C$, and $A$. Furthermore, $\tilde{d}_\Gamma(s,z)$ is a finite sum of terms of the form
$$az^b(e^{-cs})$$
for $a,b\in\intgr$ and $c\in\reals$, and $a,b=O_{g,n}(1)$. If all edge lengths of $\Gamma$ lie in $\intgr$ or $\rtnl$, then all $c$ will respectively lie in $\intgr$ or $\rtnl$ as well.
	
\end{theorem}

\begin{remark}

This theorem and the corollary below are analogous to results obtained in a more general context by Li-Matheus-Pan-Tao in their recent independent work \cite{LMPT24}, {although the hypotheses and framework used in our work are somewhat different from theirs. While it is possible the two setups can be reconciled, the author does not know of a way to do so.} The author learned about this other work while this paper was in preparation. See Subsection \ref{background} for more details.
	
\end{remark}

\subsection{Resonance Convergence and Resonance Chains}

Any approximation result on $d_X(s,z)$ and $\tilde{d}_\Gamma(s,z)$ as above can be used to obtain resonance-counting results via the argument principle. To set up this corollary, note that each $MRG_{g,n}(L_1,\dots,L_n)$ lies inside a larger moduli space $MRG_{g,n}$ where boundaries of any lengths are allowed, and in fact the individual $MRG_{g,n}(L_1,\dots,L_n)$ stratify this space.

\begin{corollary}
\label{maincorollary}

Fix $\Gamma\in MRG_{g,n}(L_1,\dots,L_n)$, and assume $\Gamma$ is trivalent. Let $(X_{{t}})$ be a sequence of surfaces with topological type $(g,n)$ (for $t\in\ntrl$), and let $(\alpha_{{t}})$ be a sequence of positive real numbers, such that the following properties hold:
\begin{enumerate}
\item $\lim_{{t}\to\infty}\alpha_{{t}}=\infty$.
\item If the sequence $(\Gamma_{{t}})$ in $MRG_{g,n}$ is defined by $\Gamma_{{t}}=\frac{1}{\alpha_{{t}}}\Phi(X_{{t}})$, then 
$$\lim_{{t}\to\infty}\frac{1}{\alpha_{{t}}}\Gamma_{{t}}=\Gamma.$$
\end{enumerate}
Also choose some connected bounded open domain $U\subset\complex$ such that $\partial U\cap \tilde{\Res}(\Gamma)=\emptyset$. Then
\begin{equation}\lim_{{t}\to\infty}|U\cap\Res(X_{{t}},\alpha_{{t}})|=|U\cap\tilde{\Res}(\Gamma)|.\end{equation}
	
\end{corollary}

Since $U$ can be arbitrarily small, this result provides a sense in which the normalized resonances of the {surfaces} $X_{{t}}$ converge to the resonances of $\Gamma$. One important consequence of such a result is that it shows that many hyperbolic surfaces with large funnel widths possess approximate \emph{resonance chains}, or sets of resonances arranged with equal spacing along a vertical line in $\complex$. {Resonance chains have been extensively studied in connection to surface Laplacians and other dynamical objects; see Subsection \ref{background} for a very brief background.}

To explain the appearance of approximate resonance chains, consider $X\in\M_{g,n}(L_1,\dots,L_n)$ and $\Gamma=\Phi(X)$. Recall the structure of $\tilde{d}_\Gamma(s,z)$ presented in Theorem \ref{mainresult}. If $\Gamma$ has all integer edge lengths, then $\tilde{d}_\Gamma(s,1)$ is a polynomial in $e^{-s}$, and in particular its root set is the union of finitely many (exact) resonance chains of the form
$$\{s_0+2\pi ik:k\in\intgr\}.$$
If instead all edges of $\Gamma$ have rational lengths, we may obtain an analogous statement by finding a common denominator and replacing $k\in\intgr$ above with $k\in\alpha\intgr$ for some well-chosen $\alpha$. By the resonance scaling property, resonance chains also exist for $\tilde{d}_{\Gamma}(s,1)$ whenever $\Gamma$ can be rescaled so that all edge lengths lie in $\intgr$ or $\rtnl$.

\subsection{Symmetric Three-Funneled Spheres}
\label{threesphere}

As an example, consider the moduli space $\M_{0,3}(L,L,L)$; this space consists of a single three-funneled sphere $X_L$. Define $\Gamma_L=\Phi(X_L)$. One may calculate that $\Gamma_L$ is a theta graph with three edges of length $\frac{L}{2}$. In particular,
$\Gamma_L=L\Gamma_1$. Furthermore, we have (see Example \ref{iharaex} for a similar calculation) that
\begin{equation}\tilde{d}_{\Gamma_1}(s,z)=-4z^6e^{-3s} + 9z^4e^{-2s} - 6z^2e^{-s} + 1\end{equation}
{Setting $z=1$, this expression factors to}
{$$-(e^{-s}-1)^2(4e^{-s}-1)$$}
and so
$$\tilde{\Res}(\Gamma_1)=\{\ln(4)+2\pi i k,2\pi i k:k\in\intgr\}.$$
In particular, the critical exponent of $\Gamma_1$ (defined here as the real part of the rightmost resonance(s)) is $\ln(4)\approx 1.386$. 

The surface $X_L$ is $\frac{L}{4}$-great by symmetry, so for any fixed $A$ and any $s,z\in\complex$ satisfying, say, $|z|\leq 2$ and $|s|\leq \frac{2}{L}$, we have {by Theorem \ref{mainresult} that}
$$|d_{X_L}(s,z)-\tilde{d}_{\Gamma_L}(s,z)|=O_\epsilon(e^{-{(1-\epsilon)L/4}}).$$
Recall that the \emph{critical exponent} $\delta_{X_L}$ of $X_L$ is likewise equal to the real part of the rightmost zero(s) of $d_{X_L}(s,z)$. By a straightforward application of Rouch{\'e}'s theorem, we obtain the following result:

\begin{corollary}
\label{funnelsphere}

Define $X_L$ to be the sole surface in the moduli space $\M_{0,3}(L,L,L)$, and let $\delta_{X_L}$ be its critical exponent. Then
\begin{equation}\delta_{X_L}=\frac{\ln(4)}{L}+O_\epsilon(e^{-{(1-\epsilon)L/4}}).\end{equation}

\end{corollary}

Versions of this result have been known for some time. Specifically, convergence was proven without speed by McMullen \cite{CM98}, was proven with error $O(L^{-3/2})$ by Pollicott-Vytnova \cite{PV19}. {Recently, a matching error of $O_\epsilon(e^{-(1-\epsilon)L/4})$} was proven by Li-Matheus-Pan-Tao \cite{LMPT24}.

\subsection{Background}
\label{background}

This work is largely based on methods developed by Weich to study resonances for three-funneled spheres \cite{TW15}. Specifically, Weich proved versions of Theorem \ref{mainresult} and Corollary \ref{maincorollary} for these surfaces, albeit without the spine graph terminology. It is primarily the relation to spine graphs that allows his methods to be extended here. While the ideas of dynamical zeta functions for holomorphic iterated function schemes were first developed by Bowen-Series and Ruelle, their application in this context is due to Jenkinson-Pollicott \cite{JP02}.

It should be noted that one of the primary motivations of analyzing surface resonances is to obtain precise estimates on the \emph{critical exponents} of these surfaces; the critical exponent is a numerical invariant controlling both the growth rate of the geodesic counting function and the Hausdorff dimension of the boundary of the universal cover. The critical exponent of a surface is equal to the value of the rightmost resonance in $\Res(X)$; which is always real \cite{SP76,SP88}. This line of inquiry was carried out by McMullen \cite{CM98}, Weich \cite{TW15}, Pollicott-Vytnova \cite{PV19}, Dang-Mehmeti \cite{DM24}, and Li-Matheus-Pan-Tao \cite{LMPT24}.

The other motivation is to explain the appearance of approximate resonance chains in surface resonance sets. Resonance chains were first observed by Borthwick \cite{DB14}, and have been studied by Weich \cite{TW15}, Li-Matheus-Pan-Tao \cite{LMPT24}, and others, as well as by Barkhofen-Kuhl-Poli-Schomerus-Weich in the context of theoretical physics \cite{BKPSW14}.

As advertised above, this work overlaps conceptually with recent independent work of Li-Matheus-Pan-Tao, trading a more limited scope for an approach some readers may find more geometrically concrete. The primary difference is in the geometric perspective taken and the model used. Specifically, Li-Matheus-Pan-Tao consider families of Schottky groups over $\complex$ degenerating to a non-Archimedean Schottky group, which respectively play the role of surfaces and spine graphs in our setup. This family is captured via a meromorphic function $z\to\Gamma_z$ into Schottky space satisfying certain specific conditions, and the limit is taken as $z\to 0$. The implicit constants in their results depend on the family $\Gamma_z$ chosen. 

While it would be very interesting to compare their bounds to ours in cases more general than that given in Subsection \ref{threesphere}, doing so seems difficult. Specifically, for a degenerating family $\Gamma_z$ consisting entirely of Schottky groups over $SL_2(\reals)$ (and thus corresponding to surfaces), one would need a way to calculate our limit parameters $L$ and $\eta$ from their limit parameter $z$. At present, we do not know of a way to do so.

\subsection{Organization}

Section two contains a more thorough overview on resonances, resolvents, the Bowditch-Epstein spine construction, and zeta functions. In section three, we define dynamical zeta functions for surfaces in terms of holomorphic iterated function schemes (IFSes). In {section} four, we study how the geometry of $\eta$-great surfaces gives some control over the data defining these IFSes. In sections five and six, we use this control along with a cycle expansion of dynamical zeta functions to establish important estimates on these functions. Lastly, in {section} seven we briefly assemble these estimates to prove our main results.

\subsection{Asymptotic Notation and Implicit Constants}

Throughout this paper, we will generally consider a surface $X$ lying in a moduli space $\M_{g,n}(L_1,\dots,L_n)$, where the $L_h$ are $C$-bounded by a parameter $L$. {The variable $\Gamma$ will almost always represent the corresponding metric ribbon graph 
$$\Gamma=\Phi(X)\in MRG_{g,n}(L_1,...,L_n)$$}
We will also consider dynamical zeta functions $d_X(s,z)$ and $\tilde{d}_\Gamma(s,z)$, and will frequently assume $|z|\leq A$ and $|s|\leq \frac{A}{L}$. There will also be a parameter $\epsilon>0$ introduced in Lemma \ref{inversecontrol}.

{Furthermore, we will fix parameters $d$ and $v$ to respectively represent the number of edges and vertices of $\Gamma$. When $\Gamma$ is trivalent, which will generally be assumed to be the case, a straightforward Euler characteristic computation shows that
$$d=6g-6+3n,\ v= 4g-4+2n$$
In particular, both $d$ and $v$ are $O_{g,n}(1)$.}

{Due to the number of constants in simultaneous use, we will often also use the Landau notations $O(\cdot)$ and $\Omega(\cdot)$. Recall that for real-valued functions $f$ and $g$ with the same domain, $f=O(g)$ is shorthand for the statement that there exists some constant $c \geq 0$ such that $f(x)\leq cg(x)$ for all $x$ in the common domain of $f$ and $g$. The notation $\Omega(\cdot)$ is shorthand for the statement that there exists some constant $c'\geq 0$ such that $f(x)\geq c g(x)$ for all $x$. When variables are subscripted to $O$ or $\Omega$, the implicit constants are allowed to depend on those variables. For example, $f=O_{g,n,C,A}(1)$ means that $f$ is uniformly bounded by some constant depending only on $g$, $n$, $C$, and $A$. Note that while big-$O$ notation is fairly standard, the definition of big-$\Omega$ notation used here is specifically the definition due to Knuth, as opposed to the alternate definition due to Hardy-Littlewood and Landau.}

\textbf{{Since almost all implicit constants in this paper will depend on $g$, $n$, $C$, $A$, and/or {$\epsilon$},} whenever big-$O$ and big-$\Omega$ notation is used, all such implicit constants will be allowed to depend on $g$, $n$, $C$, $A$, {and $\epsilon$} unless otherwise stated.} {Even with this standing assumption, these subscripts will sometimes be included anyways for clarity.}

Furthermore, whenever we work with a moduli space $\M_{g,n}(L_1,\dots,L_n)$ with the $L_h$ $C$-bounded by $L$ and an $\eta$-great $X\in\M_{g,n}(L_1,\dots,L_n)$, we will assume $\eta$ is {positive and bounded above by} $L^{{\upsilon}}$ for some $0{<\upsilon}<1$. \textbf{We will also assume that $\eta$ is sufficiently large in terms of $g$, $n$, $C$, $A$, {and $\epsilon$}.}

\subsection{Acknowledgments}

The author would like to thank Alex Wright, Richard Canary, Tobias Weich, and Wenyu Pan for helpful conversations. The author would also like to thank the anonymous referees for their feedback. Lastly, the author would also like to thank Daniel Cross and Charlotte Goist for assistance with the figures. This material is based upon work supported by the National Science Foundation Graduate Research Fellowship under Grant. No. DGE 2241144.

\section{Preliminaries}

\subsection{Zeta Functions}
\label{zetafunctions}

The spectrum of a surface $X$ is closely connected to its \emph{Selberg zeta function}.

\begin{definition}

Let $X$ be an infinite-area hyperbolic surface, and let $P_X$ be the set of all oriented closed primitive geodesics on $X$. Then the Selberg zeta function of $X$ is defined as
\begin{equation}Z_X(s)=\prod_{\gamma\in P_X}\prod_{k\geq 0}\left(1-e^{-(s+k)l(\gamma)}\right).\end{equation}
	
\end{definition}

\begin{theorem}[{Epstein-}Patterson-Perry, Theorem 1.1 \cite{EPP01}]
\label{selbergconnection}

For an infinite-area hyperbolic surface $X$, the zero set of the zeta function $Z_X(s)$ is the union of $\Res(X)$ (with matching multiplicity) and the negative integers.
	
\end{theorem}

Defining the dynamical zeta function $d_X(s,z)$ used in Theorem \ref{mainresult} requires significant additional setup, and is delayed until Lemma \ref{selbergequality}.

We define dynamical zeta functions for metric graphs as follows. Let $\Gamma$ be a finite (undirected) metric graph with {$d$} edges. Take two copies of each edge with opposite orientations and enumerate them as $\ep_1,\dots,\ep_{2 {d}}$. Say $\ep_{{j}}$ `feeds into' $\ep_{{k}}$ if the terminal vertex of $\ep_{{j}}$ is also the initial vertex of $\ep_{{k}}$, and $\ep_{{j}}$ and $\ep_{{k}}$ are not derived from the same undirected edge. 

Next, define a family of $2{d} \times 2{d}$ matrices $J_s$ parameterized by $s\in\complex$. We think of these matrices as operators on the space of functions $h:\{\epsilon_1,\dots,\epsilon_{2{d}}\}\to\complex$, defined by
\begin{equation}(J_sh)(\epsilon_{{j}})=\sum_{{k}:{j}\to {k}}e^{-s(l(\ep_{{j}})+l(\ep_{{k}}))/2}h(\epsilon_{{k}}).\end{equation}
{In matrix notation, this definition becomes}
{\begin{equation}\label{graphmatrixdef}(J_s)_{jk}=\begin{cases}e^{-s(l(\epsilon_j)+l(\epsilon_k))} & j\to k\\ 0&\text{otherwise}\end{cases}\end{equation}}
The dynamical zeta function of $\Gamma$ is then defined as
\begin{equation}\label{graphzetadef}\tilde{d}_\Gamma(s,z)=\det(I-zJ_s).\end{equation}
One can check that if all edge lengths of $\Gamma$ lie in $\intgr$, then in fact all exponents in $\tilde{d}_\Gamma(s,z)$ are of the form $e^{-cs}$ for $c\in\intgr$, rather than $c\in\frac{1}{2}\intgr$ as one might expect.

\begin{figure}[h]

\centering
\includegraphics[scale=0.8]{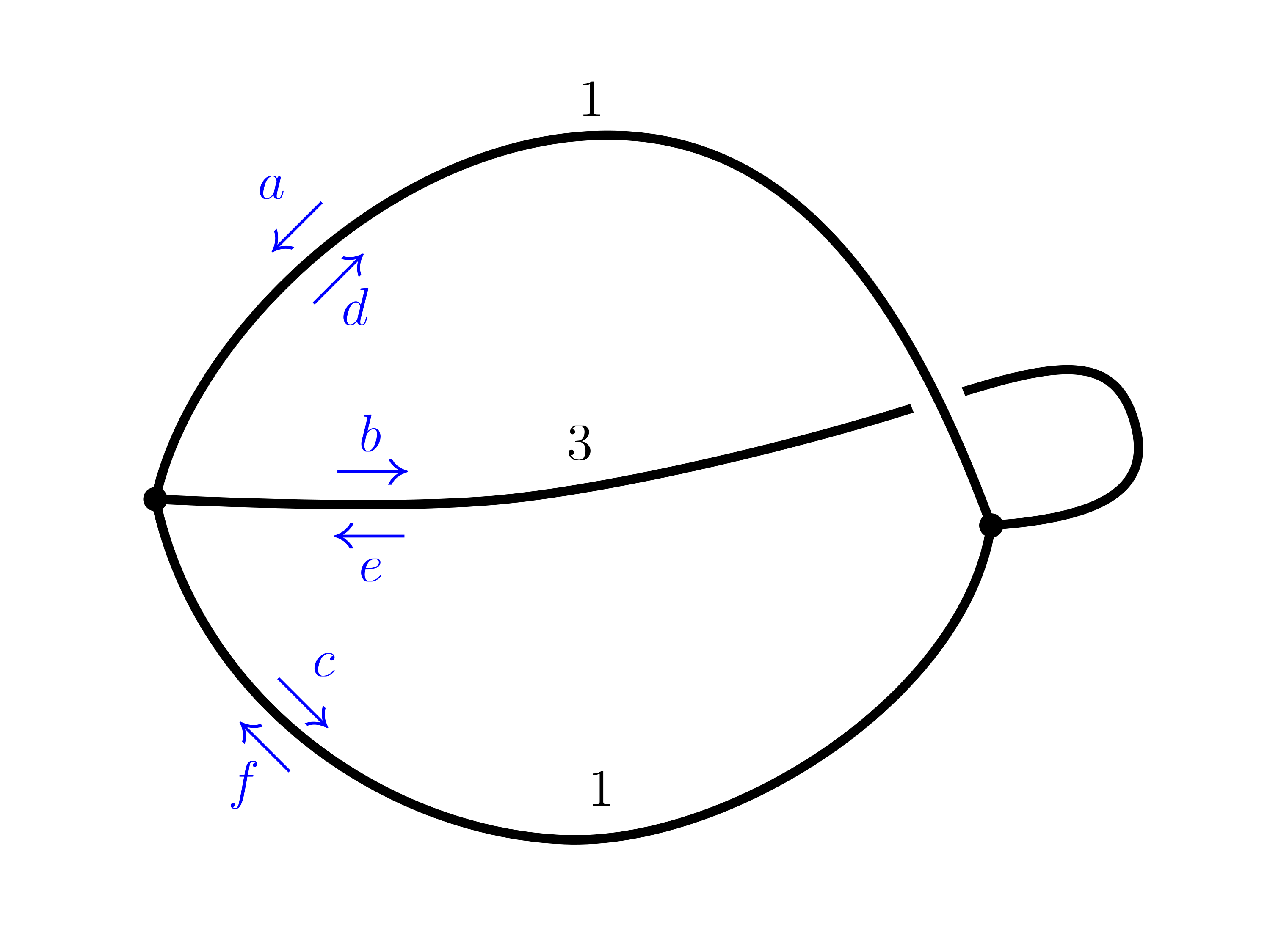}
\caption{The metric ribbon graph used in Example \ref{iharaex}, with labeled directed edges.}

\label{graphfig}

\end{figure}

\begin{example} 

\label{iharaex}

Consider the metric ribbon graph $\Gamma$ depicted in Figure \ref{graphfig}, with edge lengths $1$, $3$, and $1$. Label the directed edges as $a$, $b$, $c$, $d$, $e$, and $f$ as in Figure \ref{graphfig}. Then the matrix family $J_s$ can be written as

$$J_s=\begin{bmatrix}0&0&0&0&e^{-2s}&e^{-s}\\ 0&0&0&e^{-2s}&0&e^{-2s}\\0&0&0&e^{-s}&e^{-2s}&0\\0&e^{-2s}&e^{-s}&0&0&0\\e^{-2s}&0&e^{-2s}&0&0&0\\e^{-s}&e^{-2s}&0&0&0&0\end{bmatrix}$$
so that
\begin{equation}\tilde{d}_\Gamma(s,z)=-4z^6e^{-10s} + z^4e^{-4s} + 4z^4e^{-6s} + 4z^4e^{-8s} - 2z^2e^{-2s} - 4z^2e^{-4s} + 1.\end{equation}
Setting $z=1$, we obtain the following {polynomial:}
{$$-4(e^{-s})^{10}+4(e^{-s})^8+4(e^{-s})^{6}-3(e^{-s})^{4}-2(e^{-s})^2+1$$}
{$$=-(e^{-s}-1)^2(e^{-s}+1)^2(4(e^{-s})^6+4(e^{-s})^4-1)$$}
{We then obtain the following} resonance chains for $k\in\intgr$:
$$s=\pi i k,\ s\approx 0.434+\pi i k,\ s\approx 0.129\pm 1.369i+\pi i k.$$
In particular, the critical exponent of $\Gamma$ is approximately $0.434$.

\end{example}

\subsection{Spines}

For any $X\in \M_{g,n}(L_1,\dots,L_n)$, let $X_C$ be its compact core; i.e. the region lying between all funnel geodesics. Note that any closed geodesic on $X$ must lie entirely in $X_C$. Define the spine of $X$ to be the set of points $x\in X_C$ such that $d(x,\partial X_C)$ is realized by multiple points on $\partial X_C$. The spine has a nice geometric structure.

\begin{theorem}[Bowditch-Epstein \cite{BE88}]
\label{spineconst}

Fix $X\in \M_{g,n}(L_1,\dots,L_n)$, let $E$ be the set of $x\in X_C$ such that $d(x,\partial X_C)$ is realized by exactly two boundary points, and let $V$ be the set of $x\in X_C$ such that $d(x,\partial X_C)$ is realized by at least three boundary points. Then $E\cup V$ is an embedded metric ribbon graph, the spine graph $S_X$ of $X$, with $V$ as its vertex set. Furthermore, $X_C$ deformation retracts onto $S_X$.
	
\end{theorem}

Here, a \emph{metric ribbon graph} is a finite metric graph with an orientation around each vertex; one thinks of inflating each edge into a thin ribbon and using the vertex orientations to glue the ribbons together into a thin (topological) surface. {For a metric graph embedded in an (oriented) surface, the orientation naturally induces a ribbon structure on the graph.}

\begin{definition}

The space $MRG_{g,n}(L_1,\dots,L_n)$ is the moduli space of metric ribbon graphs with genus $g$, $n$ boundary components, and prescribed boundary lengths $L_1,\dots,L_n$. Here, metric ribbon graphs are identified up to metric ribbon graph isometry.
	
\end{definition}

\begin{remark}

Like $\M_{g,n}(L_1,\dots,L_n)$, the space $MRG_{g,n}(L_1,\dots,L_n)$ possesses a rich geometry. For more information, see \cite{ABCGLW21}.
	
\end{remark}

We can now construct a surface-to-spine map
\begin{equation}\Phi:\M_{g,n}(L_1,\dots,L_n)\to MRG_{g,n}(L_1,\dots,L_n).\end{equation}
Topologically, $\Phi(X)$ will be the same as the spine graph $S_X$ defined above. However, defining the lengths of its edges is more subtle. First draw for each $x\in V$ and $y\in \partial X_C$ realizing $d(x,\partial X_C)$ the geodesic arc from $x$ to $y$. These arcs are the \emph{ribs} of the spine, and divide $X_C$ into a union of hexagons we will call \emph{corridors}. Each corridor contains one edge of $S_X$ and two boundary arcs of $\partial X$ with equal length; we define the length of the corresponding edge in $\Phi(X)$ to be the length of either boundary arc. The key property of this construction is that the lengths of the boundary faces of $\Phi(X)$ are identical to the funnel widths of $X$, so $\Phi$ has the correct codomain.

\begin{theorem}[Bowditch-Epstein, Theorem 9.5 \cite{BE88}]
\label{spinemap}

The function $\Phi:\M_{g,n}(L_1,\dots,L_n)\to MRG_{g,n}(L_1,\dots,L_n)$ defined above is a homeomorphism.
	
\end{theorem}

It will be handy to decompose {a surface core $X_C$} into {polygons} in a slightly different way: for each corridor, one may draw an arc contained in the corridor between the two boundary components of that corridor, and then find the minimal-length representative of the homotopy class of that arc relative to the boundary. These representatives are called \emph{intercostals}, and cut a given surface core into {polygons} we will call \emph{sectors}. If the intercostals are extended into bi-infinite rays on $X$, they instead divide $X$ into $\emph{long sectors}$. Each long sector is isometric to a region in $\H$ bounded between some number non-intersecting geodesics. {When the spine graph $S_X$ of $X$ is trivalent, which will almost always be the case in our analysis, these sectors will be hexagons, and long sectors will be bounded by three geodesics.}

\begin{figure}[h]

\centering
\includegraphics[scale=0.18]{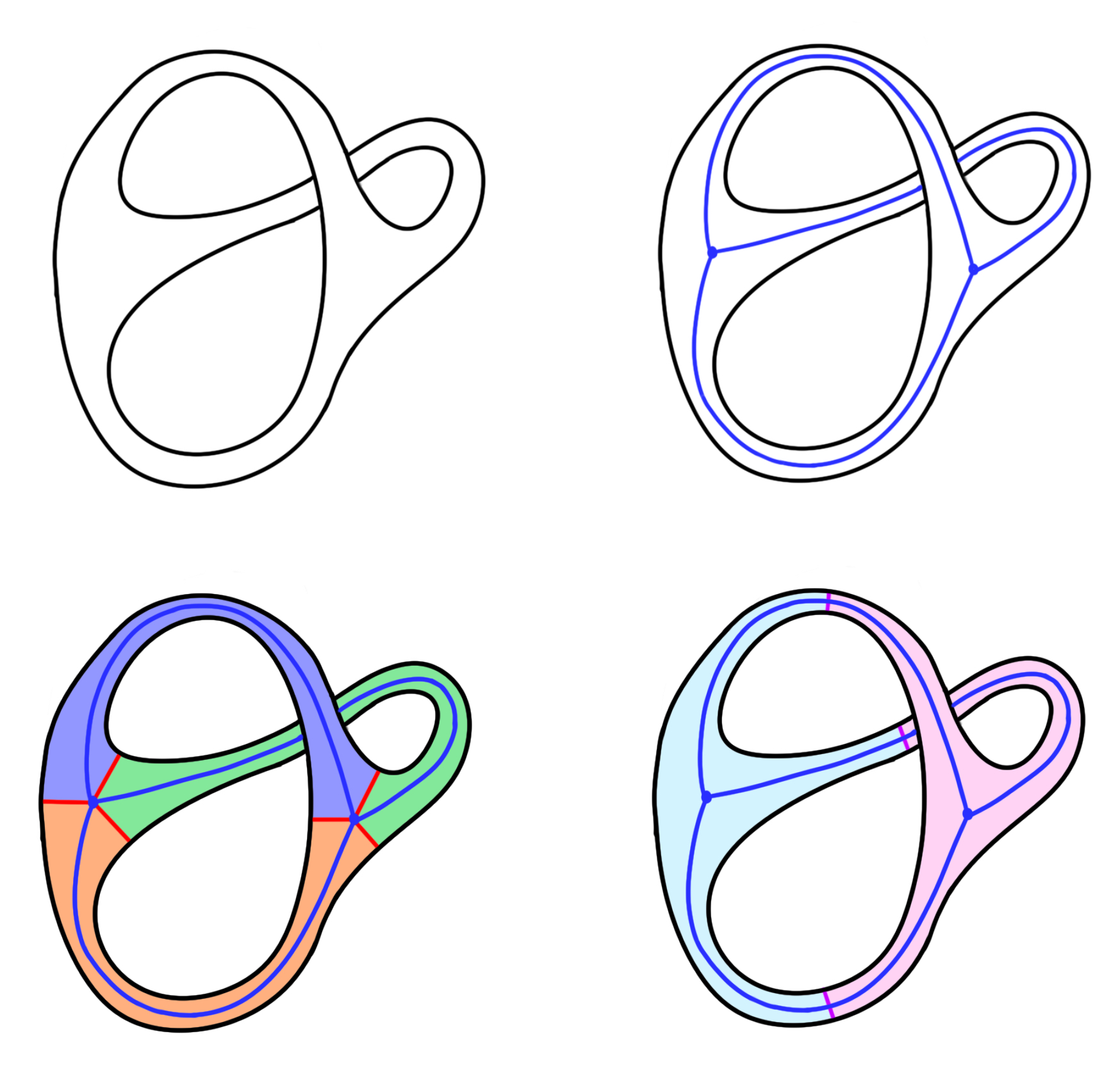}
\caption{Examples of constructions used in this paper. Top left: a surface core $X_C$ with $g=1$ and $n=1$ (i.e., a one-holed torus). Top right: the same surface core with its spine graph marked. Notice that this metric ribbon graph has a single face. Bottom left: the same surface core, with its ribs added in red and the three corridors shaded. Bottom right: the same surface core, with its intercostals added in magenta and its two sectors shaded.}

\label{spinefig}

\end{figure}

For an illustration of the constructions made in this subsection, see Figure \ref{spinefig}.

\section{Iterated Function Schemes and Dynamical Zeta Functions}

\subsection{Holomorphic Iterated Function Schemes}

We would like to reformulate resonances in terms of a dynamical system. Our primary tools for doing so will be holomorphic iterated function schemes, transfer operators, and dynamical zeta functions. While we will use a slightly modified definition of holomorphic iterated function schemes, we first present the standard definition:

\begin{definition}

A \emph{holomorphic iterated function scheme} (holomorphic IFS) consists of the following data:
\begin{enumerate}
\item A collection of open disks $D_1,\dots,D_N\subset\complex$ with pairwise disjoint closures.
\item A $N\times N$ matrix $A$ with all entries equal to $0$ or $1$. We denote $A_{{jk}}=1$ by ${j}\to {k}$.
\item For every $({j},{k})$ pair with ${j}\to {k}$, a biholomorphism $f_{{jk}}:D_{{j}}\to D_{{k}}$, such that
$$f_{{jk}}(D_{{j}})\cap f_{{j'k'}}(D_{{j'}})=\emptyset$$
unless $({j},{k})=({j'},{k'})$. We also define
$$D_{{jk}}=f_{{jk}}(D_{{j}}).$$
\end{enumerate}
	
\end{definition}

We will modify this definition in two ways. First, we will introduce spaces $C_1,\dots,C_N$, each a copy of $\complex$, and think of $D_{{j}}$ as lying in $C_{{j}}$. The first condition above no longer makes sense, so may be discarded. Second, we will replace the third condition with the condition that for all ${j},{j'}$ and all ${k}$ with ${j}\to {k}$, ${j'}\to {k}$,
$$f_{{jk}}(D_{{j}})\cap f_{{j'k}}(D_{{j'}})=\emptyset$$
unless ${j}={j'}$.

For us, holomorphic IFSes will be the only type of IFSes under consideration, so we will sometimes drop the `holomorphic' qualifier.

\begin{remark}

Note that any such IFS can be modified into one that matches the original definition by translating each $D_{{j}}$ so that they have disjoint closures when projected to a common copy of $\complex$, performing the projection, and then replacing the $f_{{jk}}$ with $\tilde{f}_{{jk}}=f_{{jk}}+z_{{jk}}$ for some carefully chosen translation terms $z_{{jk}}$. Modifying the $f_{{jk}}$ in this way has no effect on the IFS properties we care about (in particular, the derivatives of the $f_{{jk}}$ are unchanged).

\end{remark}

\begin{remark}

While it is not strictly part of the definition, all IFSes we will consider will be such that the disks have centers in $\reals$, and all $f_{{ij}}$ will be linear fractional transformations with real coefficients.

\end{remark}

An IFS generates a natural symbolic coding. The symbols of this coding are integers $1,\dots,N$, and a word of length ${m}$ in this coding is a sequence $w_0\cdots w_{{m}}$ such that $w_{{m'}}\to w_{{m'}+1}$ for all $0\leq {m'}\leq {m}-1$. Each word corresponds to a map
$$f_w=f_{w_{{m}-1}w_{{m}}}\circ\dots\circ f_{w_0w_1}:D_{w_0}\to D_{w_{{m}}}.$$
We define $W_{{m}}$ to be the set of all words of length ${m}$, and $W_{{m}}^C$ to be the set of words of length ${m}$ satisfying $w_0=w_{{m}}$. Words satisfying this second condition are referred to as \emph{closed words}.

For a given IFS {with disks $D_1,...,D_N$}, we will also define
$$D=\bigcup_{{j}=1}^{{N}} D_{{j}}$$
and
$$f(D)=\bigcup_{{j}\to {k}}D_{{jk}}.$$
Notice that $f(D)$ is a disjoint union under our definition. 

Also, by the disjoint images condition and biholomorphicity, we may define a partial inverse map
\begin{equation}\label{bowenserieseq}f^{-1}:f(D)\to D,\ f^{-1}|_{D_{{jk}}}=f_{{jk}}^{-1}.\end{equation}

\subsection{Transfer Operators}

Given an IFS with disks $D_1,\dots,D_N$, let $A_\infty(D)$ be the Banach space of holomorphic functions on $D$ that extend to continuous functions on the closure $\overline{D}$ and are bounded on $\overline{D}$ with the supremum norm.

\begin{definition}

Let $V\in A_\infty(f(D))$ be a chosen \emph{potential function}. Then the \emph{transfer operator} ${\L}_V:A_\infty(D)\to A_\infty(D)$ associated to the given IFS is defined by
$${\L}_Vh|_{D_{{j}}}={\L}_V^{{j}}h,$$
$${\L}^{{j}}_V:A_\infty(D)\to A_\infty(D_{{j}}),$$
\begin{equation}({\L}^{{j}}_Vh)(u)=\sum_{{k}:{j}\to {k}} V(f_{{jk}}(u))h(f_{{jk}}(u)).\end{equation}
	
\end{definition}

We also define iterated potentials corresponding to words. For $w\in W_{{m}}$,
\begin{equation}V_w(u)=\prod_{{m'}=1}^{{m}}V(f_{w_{0,{m'}}}(u))\end{equation}
where $w_{0,{m'}}$ is the truncated word obtained from the first ${m'}+1$ letters of $w$.

With this definition, a straightforward calculation shows that
\begin{equation}({\L}_V^{{m}}h)(u)=\sum_{w\in W_{{m}},\ u\in D_{w_0}} V_w(u)h(f_w(u)).\end{equation}

\begin{definition}

An operator ${\L}:B\to B$ on a Banach space $B$ is a \emph{nuclear operator} if there exist sequences $v_{{m}}\in B$, $\alpha_{{m}}\in B^*$, and $\lambda_{{m}}\in\complex$ such that $\|v_{{m}}\|=\|\alpha_{{m}}\|=1$ for all ${m}$, $\sum_{{m}=0}^\infty|\lambda_{{m}}|<\infty$, and ${\L}$ has the series representation
\begin{equation}{\L}h=\sum_{{m}=0}^\infty \lambda_{{m}}\alpha_{{m}}(h)v_{{m}}.\end{equation}
	
\end{definition}

\begin{definition}
A holomorphic IFS is \emph{eventually contracting} if there exists some $\theta<1$ such that for all sufficiently long words $w$,
$$|f_w'(u)|\leq\theta.$$
\end{definition}

Note that if a holomoprhic IFS is eventually contracting, then by taking powers, any closed word $w$ is such that the function $f_w:D_{w_0}\to D_{w_0}$ must have a unique fixed point in its domain. We denote this fixed point as $u_w$.

\begin{lemma}[Jenkinson-Pollicott, Proposition 2 \cite{JP02}]

Consider an eventually contracting holomorphic IFS. The corresponding transfer operator (for any potential function $V\in A_\infty(D)$ as defined above) is nuclear.

\end{lemma}

\begin{definition}

We define the \emph{dynamical zeta function} associated to a nuclear transfer operator $\L_V$ via the Fredholm determinant:
\begin{equation}d_V(z)=\det(1-z{\L}_V).\end{equation}
	
\end{definition}

These dynamical zeta functions satisfy a trace formula:

\begin{lemma}[Jenkinson-Pollicott, Equation 3.26 \cite{JP02}]
\label{expform}

Let $d_V(z)$ be the dynamical zeta function associated to an eventually contracting holomorphic IFS. Then $d_V(z)$ is an entire function, and the following formula holds for $|z|$ sufficiently small:
\begin{equation}d_V(z)=\mathrm{exp}\left(-\sum_{{m}>0}\frac{z^{{m}}}{{m}}\sum_{w\in W^C_{{m}}}V_w(u_w)\frac{1}{1-f'_w(u_w)}\right)\end{equation}
where $W^C_{{m}}$ is the set of closed words of length ${m}$ and $u_w$ is the unique fixed point of a given $f_w$.

\end{lemma}

\subsection{Geodesic Coding}

Fix $X\in\M_{g,n}(L_1,\dots,L_n)$. We want to set up a holomorphic IFS corresponding to $X$ in a particular way, so that it aligns with the spine graph $\Gamma=\Phi(X)$. Assume that $\Gamma$ is trivalent. Since $X$ deformation retracts onto its spine $S_X$ and $S_X$ is naturally homeomorphic to $\Gamma=\Phi(X)$, there is a bijection between (directed) geodesics on $X$ and (directed) geodesics on $\Gamma$. Now, take the edges $e_1,\dots,e_{{d}}$ of $\Gamma$. We will arbitrarily orient each edge, and use $e_{{j}}$ for $1\leq {j}\leq {d}$ to refer to the oriented edge and $e_{{j}+{d}}$ to refer to that edge with the opposite orientation. This gives a collection $e_1,\dots,e_{2{d}}$ of all oriented edges on $X$. In general, for an edge $e_{{j}}$ with $1\leq {j}\leq 2{d}$, we will use the notation $e_{{j}}^{-1}$ for that edge with the opposite orientation. 

Define a $2{d}\times 2{d}$ matrix $A$ by setting $A_{{jk}}=1$ if $e_{{j}}$ feeds into $e_{{k}}$ in the sense of Subsection \ref{zetafunctions}, and $A_{{jk}}=0$ otherwise. A length-${m}$ word $w=w_0\cdots w_{{m}}$ in the alphabet $1,\dots,2{d}$ is \emph{permissible} if for every consecutive ${jk}$ pair in $w$, $A_{{jk}}=1$. Recall that a word $w=w_0\cdots w_{{m}}$ is \emph{closed} if $w_0=w_{{m}}$.

Define the combinatorial length of a closed geodesic on $\Gamma$ as the number of edges it traverses, counting with multiplicity. Then there is a clear map from closed permissible words of length ${m}$ {to} closed geodesics on $\Gamma$ of combinatorial length ${m}$ by following directed edges, and this map is a bijection if one replaces closed permissible words with cyclic permutation classes of closed permissible words. If $w=w_0\cdots w_{{m}}$, then we use $\tilde{\gamma}_w$ to denote the corresponding geodesic:
$$[w_0\cdots w_{{m}}]\leftrightarrow\tilde{\gamma}_w.$$

Of course, since there is a bijection $\gamma\leftrightarrow\tilde{\gamma}$ between geodesics on $X$ and geodesics on $\Gamma$, we also get a bijection between cyclic permutation classes of closed permissible words and closed geodesics on $X$:
$$[w_0\cdots w_{{m}}]\leftrightarrow \gamma_w.$$
We may interpret this relation as follows: each directed edge on $\Gamma$ corresponds to a directed intercostal on $X$, i.e. an intercostal together with a perpendicular vector field allowing us to distinguish one side of the intercostal. Then any oriented geodesic $\gamma$ on $X$ induces a word by tracking which intercostals it passes through and recording its direction at each such intercostal.

\emph{Powers} of closed words are obtained by concatenation, and a primitive word is one that is not a power of a shorter word.

\subsection{The Flow-Adapted Iterated Function Scheme}
\label{flowifs}

Our goal is to build a holomorphic iterated function scheme derived from a surface $X\in\M_{g,n}(L_1,\dots,L_n)$. Since this iterated function scheme will follow a specific construction (and other ways of associating a surface to an IFS have been defined in the literature), we will refer to it as a \emph{flow-adapted} iterated function scheme, following terminology used by Weich in \cite{TW15}.

Fix $X\in\M_{g,n}(L_1,\dots,L_n)$, let $S_X$ be its spine, and set $\Gamma=\Phi(X)$. Assume $\Gamma$ is trivalent, and that all corridors of $S_X$ contain their intercostals - this is true for $\eta$-great surfaces {by definition}. $X_C$ admits a decomposition into hexagonal sectors, which extend to `long sectors' on $X$ by extending intercostals. Each sector is bounded by three boundary arcs and three intercostal arcs, and the intercostal arcs on $X_C$ extend to bi-infinite geodesics we term \emph{long intercostals} on $X$.

Let $e_1,\dots,e_{{d}}$ be the edges of $\Gamma$; it will be convenient for the entirety of this section to think of $\Gamma$ as topologically embedded in $X$ via identification with $S_X$. We arbitrarily orient each edge, and define $e_{{d}+j}$ to be $e_j$ with reversed orientation for $1\leq j\leq {d}$. This gives an enumeration $e_1,\dots,e_{2{d}}$ of all possible directed edges of $\Gamma$. The disks $D_1,\dots,D_{2{d}}$ of our IFS will be in bijection with these directed edges.

Fix some $e_j$, and define $H_j\subset X$ to be the long sector in which $e_j$ originates. The edge $e_j$ intersects exactly one intercostal, which we will denote as $\iota_j$. This intercostal intersects $\partial X_C$ in two points, which we will denote as $l_j$ and $r_j$. The labels are chosen so that if the sector $H_j\cap X_C$ is oriented clockwise, $l_j$ lies to the left of $r_j$ along $\iota_j$.

\begin{figure}[h]

\centering
\includegraphics[scale=0.25]{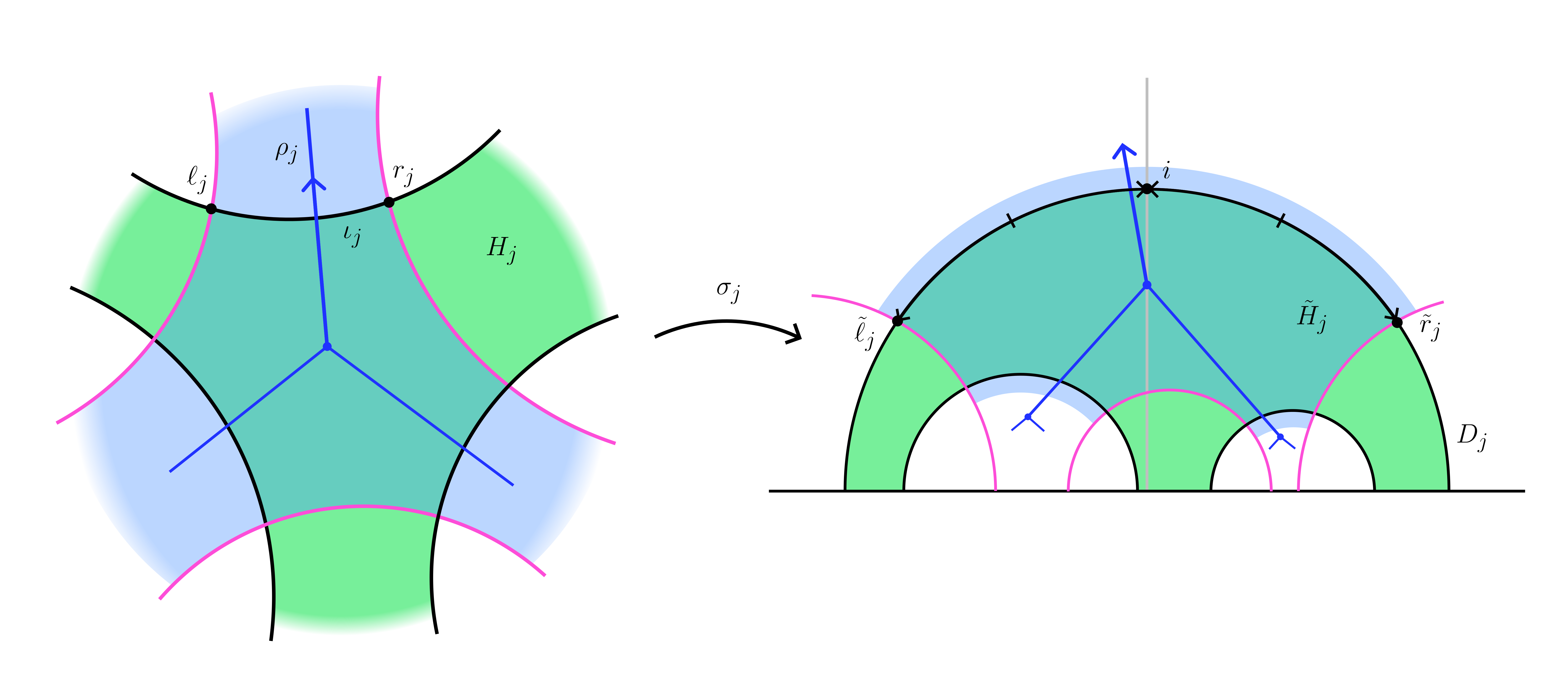}
\caption{An illustration of the map $\sigma_j$ used in the construction of the IFS. On the left, the green shading represent the long sector $H_j$, while the blue shading represents the core surface $X_C$. On the right, the colors represent the images of these regions under $\sigma_j$. Note that while the entirety of the image of $H_j$ is depicted, the entirety of the image of $X_C$ is not shown.}

\label{sigmafig}

\end{figure}

Consider a copy $\H_j$ of $\H$ defined using the upper-half-plane model. Let 
$$\sigma_j:H_j\hookrightarrow\H_j$$
be an orientation-preserving, holomorphic, isometric embedding with the following properties:
\begin{enumerate}
\item The intercostal $\iota_j$ is mapped onto a subarc of the geodesic between $-1$ and $1$ on $\H_j$ (i.e., the {Euclidean} half-circle with center $0$ and radius $1$).
\item The long sector $H_j$ is mapped into the interior of the half-disk bounded by this geodesic.
\item The points $\tilde{l}_j=\sigma(l_j)$ and $\tilde{r}_j=\sigma(r_j)$ are equidistant from {the point} $i$ along $\sigma(\iota_j)$.
\end{enumerate}
These conditions can always be satisfied, and in fact define $\sigma_j$ uniquely. The construction of the map $\sigma_j$ is demonstrated in Figure \ref{sigmafig}.

We may now embed $\H_j$ in $\complex$, and define $D_j$ to be the disk centered at $0$ of radius $1$ in this copy of $\complex$. We also define $\tilde{H}_j=\sigma_j(H_j)$, $\tilde{\iota}_j=\sigma_j(\iota_j)$, $\tilde{l}_j=\sigma_j(l_j)$, and $\tilde{r}_j=\sigma_j(r_j)$.

The transition matrix $A$ is defined as in the previous subsection, where $A_{{jk}}=1$ if $e_{{j}}$ feeds into $e_{{k}}$ and $A_{{jk}}=0$ otherwise.

Next, we need to define the transition maps $f_{{jk}}$. Letting $\tilde{X}$ be the universal cover of $X$, the first step is to note that $\sigma_{{k}}:H_{{k}}\hookrightarrow\H_{{k}}$ extends uniquely to an isometry $\sigma_{{k}}:\tilde{X}\to\H_{{k}}$. We think of these isometries as identifications and use them to define covering maps $\pi_{{k}}:\H_{{k}}\to X$.

For any ${j}\to {k}$, we observe that $e_{{j}}$ originates in $H_{{j}}$ and terminates in $H_{{k}}$, and in particular $H_{{j}}$ and $H_{{k}}$ are adjacent and separated by $\iota_{{j}}$. As a consequence, $\tilde{H}_{{j}}$ is adjacent to some lift $\tilde{H}_{{jk}}$ of $H_{{j}}$ (lifting with respect to $\pi_{{k}}$), and separated from this lift by a lift $\tilde{\iota}_{{jk}}$ of $\iota_{{j}}$. Since $\iota_{{j}}\neq\iota_{{k}}$, this $\tilde{\iota}_{{jk}}$ must be one of the two bounding geodesics of $\tilde{H}_{{k}}$ contained strictly within $D_{{k}}$. Since the long sectors $\tilde{H}_{{j}}$ and $\tilde{H}_{{jk}}$ are isometric, there is a unique hyperbolic isometry $F_{{jk}}:\tilde{H}_{{j}}\to\tilde{H}_{{jk}}$, and this isometry extends uniquely to an isometry
$$F_{{jk}}:\H_{{j}}\to\H_{{k}}$$
making the following diagram commute:
$$\begin{tikzcd}
\mathbb{H}_{{j}} \arrow[rr, "F_{{jk}}"] \arrow[rd, "\pi_{{j}}"] &     & \mathbb{H}_{{k}} \arrow[ld, "\pi_{{jk}}"'] \\
                                                      & X &                                  
\end{tikzcd}$$
As a complex map, $F_{{jk}}$ is a linear fractional transformation, and so extends to a meromorphic map from $\hat{\complex}$ to itself. Furthermore, its unique pole lies outside of $D_{{jk}}$, and if $D_{{jk}}$ is the disk in $\complex$ half-bounded by the geodesic extension of $\iota_{{jk}}$, then $F_{{jk}}(D_{{j}})=D_{{jk}}$. The second fact is clear by construction, while the first straightforwardly follows from the explicit formula proven in Lemma \ref{derivcontrol}. So we may define
$$f_{{jk}}:D_{{j}}\to D_{{k}},\ f_{{jk}}=F_{{jk}}|_{D_{{j}}}.$$

We associate two quantities to this map. The first is $\delta_{{j}}=d_\H(\tilde{l}_{{j}},\tilde{r}_{{j}})$, while the second is $\kappa_{{jk}}=d_\H({\partial} D_{{k}},{\partial} D_{{jk}})$. Both distances are realized by arcs of the boundary of the sector contained in the long sector $H_{{k}}$. Specifically, the distance $\delta_{{j}}$ is realized by the intercostal arc crossed by $e_{{j}}$, while the distance $\kappa_{{jk}}$ is realized by the boundary arc incident to the intercostal arc crossed by $e_{{j}}$ and the intercostal arc crossed by $e_{{k}}$. 

\begin{figure}[h]

\centering
\includegraphics[scale=0.3]{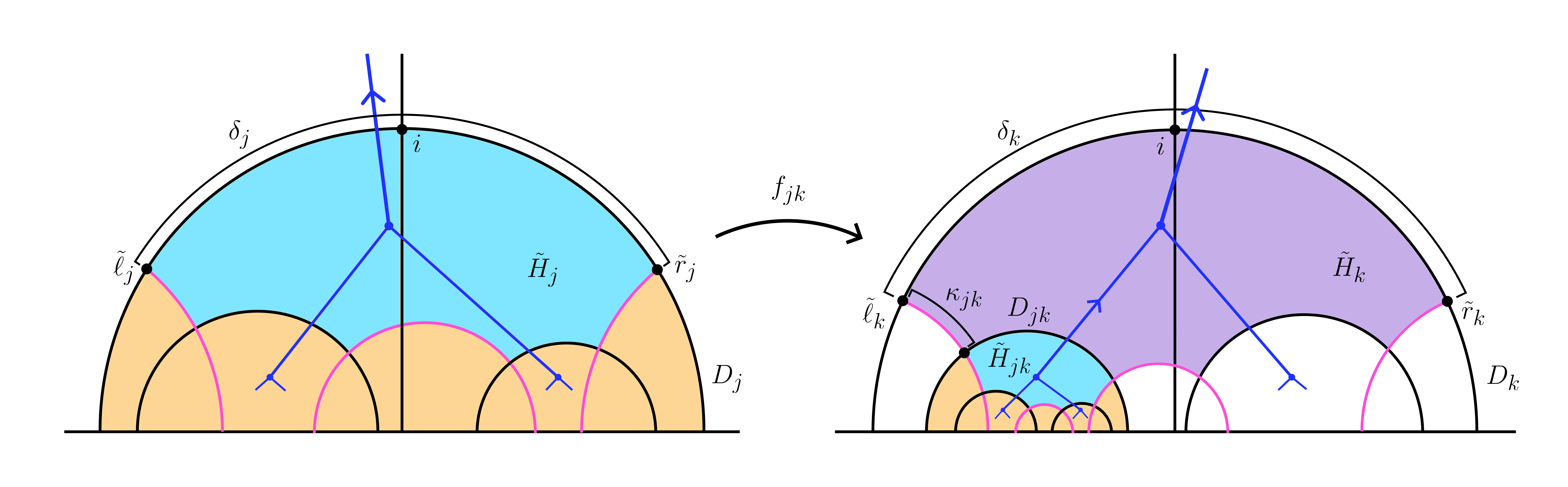}
\caption{An illustration of the map $f_{{jk}}$ and associated objects and quantities, including the quantities $\delta_{{j}}$, $\delta_{{k}}$, and $\kappa_{{jk}}$. Here, the orange shading on the left represents $D_{{j}}$, while the orange shading on the right represents its image $D_{{jk}}=f_{{jk}}(D_{{j}})$. Note that $f_{{jk}}$ is not necessarily the unique \emph{linear} map sending $D_{{j}}$ to $D_{{jk}}$; it may be a more complex isometry.}

\label{mapfig}

\end{figure}

For an illustration of the maps $f_{{jk}}$ and the definitions of the quantities $\delta_{{j}}$, $\delta_{{k}}$, and $\kappa_{{jk}}$, see Figure \ref{mapfig}.

\begin{remark}

Note that the illustration provided in Figure \ref{mapfig} is not a good representation of what happens in the long-boundary-length regime. Specifically, in our setup, we expect $\delta_{{j}}$ and $\delta_{{k}}$ to be very small and $\kappa_{{jk}}$ to be very large. See Section \ref{geometriclimits} for more information.
	
\end{remark}

For any fixed ${k}$, trivalency implies there are exactly two indices ${j}_1$ and ${j}_2$ such that ${j}_1\to {k}$ and ${j}_2\to {k}$. Furthermore, the disks $D_{{j}_1{k}}$ and $D_{{j}_2{k}}$ are (half-)bounded by distinct intercostals of $\tilde{H}_{{k}}\subset \H_{{k}}$, {as otherwise $e_{{j}_1}$ and $e_{{j}_2}$ would be equal}. As a consequence, 
$$f_{{jk}}(D_{{j}})=f_{{j'}{k'}}(D_{{j'}})\Rightarrow ({j},{k})=({j'},{k'}).$$

Each function $f_w$ is a linear fractional transformation with real coefficients, and is a biholomorphism from $D_{w_0}$ to $D_w$. While the eventually contracting property holds, its proof will be postponed until certain calculations are made in Section \ref{geometriclimits}.

How do these maps relate to geodesics? Say we have a permissible closed word $w=w_0\cdots w_{{m}}$ identified with a closed geodesic $\tilde{\gamma}_w$ on $\Gamma$; we think of $\tilde{\gamma}_w$ as lying in $S_X$. Then it is straightforward to check by stitching together commuting diagrams of the form presented above that $f_w$ maps $D_{w_0}$ to itself, and extends to a $\pi_{w_0}$-invariant isometry on $\H_{w_0}$. As a consequence, this extension of $f_w$ is a deck transformation, and corresponds up to conjugacy to a unique geodesic on $X$.

\begin{lemma}
\label{geodcorresp}

{With the setup given in the previous paragraph, the geodesic corresponding to $f_w$ is $\gamma_w$, where $\gamma_w$ is also the geodesic corresponding to $\tilde{\gamma}_w$ under the previously established correspondence between geodesics on $X$ and geodesics on $\Gamma$.}

\end{lemma}

\textbf{Proof.} {Choose $j_0,\dots,j_m,j_{m+1}$ with $j_{m+1}=j_0$ so that $w_{{q}}$ corresponds to the edge pointing from the sector $H_{{j_q}}$ to the sector $H_{{j_{q+1}}}$.} Start with the unique lift to $D_{j_0}$ of the spine vertex contained in $H_{{j}_0}$. After applying $f_{w_0w_1}$, we can follow a lift of $e_{w_0}$ from the image of the lift of this vertex to the lift of the unique spine vertex $v_{w_{{1}}}$ contained in $H_{{j}_1}$; this forms a path $\tilde{\gamma}_{{1}}$. Now, given $\tilde{\gamma}_{{q}}$ ending at a lift of the unique spine vertex contained in $H_{{j_{q}}}$, we can take its image under $f_{w_{{q}}w_{{q}+1}}$ and add a lift of $e_{w_{{q}+1}}$ to a lift of the unique spine vertex in $H_{w_{{q}+1}}$ to obtain $\tilde{\gamma}_{{q}+1}$. Iterating this procedure, we obtain $\tilde{\gamma}_{{m}+1}$, a partial lift of the closed path $\tilde{\gamma}$ on $X$, such that $\tilde{\gamma}_{{m}+1}$ originates at $f_w(v_{w_0})$ and terminates at $v_{w_0}$. {(Recall that we have embedded $\Gamma$ in $X$, so that $\tilde{\gamma}$ can be thought of as both a geodesic on $\Gamma$ and a path on $X$.)} {It is now straightforward to verify that} that the $f_w$-orbit of $v_{w_0}$ is of bounded distance from the lift of $\gamma$ passing through $\tilde{H}_{w_0}$, and so $f_w$ corresponds to $\gamma$. $\square$

\subsection{Dynamical Zeta Functions for IFSes}

Our reason for setting up this IFS is as follows. {The function $d_X(s,z)$ defined in Lemma \ref{selbergequality} for a surface $X\in\M_{g,n}(L_1,...,L_n)$ is the aforementioned \emph{dynamical zeta function} for $X$.}

\begin{lemma}
\label{selbergequality}

Fix $X\in\M_{g,n}(L_1,\dots,L_n)$, construct the flow-adapted IFS corresponding to $X$, and define the following family of potential functions:
\begin{equation}V_s(u)=((f^{-1})'(u))^{-s}.\end{equation}
Here, $f^{-1}$ is defined as in Equation \ref{bowenserieseq}. {Let $\L_{V_s}$ be the transfer operator associated to this IFS and to the potential function $V_s(u)$.} Setting $d_X(s,z)=d_{V_s}(z)$, and letting $Z_X(s)$ be the Selberg zeta function for $X$, we have that $d_X(s,z)$ is an entire function in both variables and
$$d(s,1)=Z_X(s).$$

\end{lemma}

\textbf{Proof.} This is a standard argument; see Theorem 3.5 in \cite{TW15} and Theorem 15.8 in \cite{DB16}. As such, our proof will gloss over the technical details.

The first step is the following formula (see Proposition 15.5 in \cite{DB16}): for any $w\in W^C_{{m}}$, if $u_w$ is the unique fixed point of $f_w$, then
$$f'(u_w)=e^{-l(\gamma_w)}.$$

Next, calculate that $V_{w,s}(u_w)$ and $f'_w(u_w)$ are independent of a choice of representative in $[w]$: this follows from the geometric interpretation in terms of translation distance given above. Furthermore, it is straightforward to calculate that
$$V_{w^j,s}(u_{w^j})=(V_{w,s}(u_w))^j,\ f'_{w^j}(u_{w^j})=(f'_w(u_w))^j$$
where $w^j$ is the $j$-fold concatenation of $w$.
We apply a geometric series expansion to the formula for $d(s,z)$ obtained in Lemma \ref{expform}, which holds for $|z|$ sufficiently small:
\begin{equation}d(s,z)=\exp\left(-\sum_{j\geq 0}\sum_{{m}>0}\sum_{w\in W^{C}_{{m}}}\frac{z^{{m}}}{{m}}V_{w,s}(u_w)(f'_w(u_w))^j\right).\end{equation}

We now wish to replace the sum $\sum_{{m}>0}\sum_{w\in W^{C}_{{m}}}$ over all closed words with the sum $\sum_{r>0}\sum_{[w]\in[W^{P}]}$ over all powers of equivalence classes of prime words. If a given equivalence class $[w]$ of length ${m}$ has prime representation $[w]=[w_p^r]$ where $w_p$ is prime and has length ${m'}$, then this equivalence class will be counted $r$ times under the first sum. Under the second sum, it will be counted once. Since necessarily ${m}=r{m'}$, we obtain
\begin{equation}d(s,z)=\exp\left(-\sum_{{j}\geq 0}\sum_{r>0}\sum_{[w]\in[W^{P}]}\frac{\left(z^{c(w)}V_{w,s}(u_w)(f'_w(u_w))^{{j}}\right)^r}{r}\right)\end{equation}
where $c(w)$ is the word length of $w$. We now exchange the order of summation and apply the Taylor expansion for $\log(1-x)$ to obtain
\begin{equation}d(s,z)=\prod_{[w]\in[W^{P}]}\prod_{j\geq 0}\left(1-z^{c(w)}V_{w,s}(u_w)(f'_w(u_w))^j\right).\end{equation}
We now apply the chain rule and the inverse derivative formula. {For a word $w$ with length ${m}$,}
$$V_{w,s}(u_w)=\prod_{{q}=1}^{{m}}V_s(f_{w_{0,{q}}}(u_w))$$
$$=\prod_{{q}=1}^{{m}}((f^{-1})'(f_{w_{0,{q}}}(u_w)))^{-s}=f'_w(u_w)^s.$$
So
\begin{equation}d(s,z)=\prod_{[w]\in[W^{P}]}\prod_{j\geq 0}\left(1-z^{c(w)}(f'_w(u_w))^{j+s}\right).\end{equation}
Applying the translation length formula gives
\begin{equation}d(s,z)=\prod_{\gamma}\prod_{j\geq 0}\left(1-z^{c(w)}e^{-l(\gamma)(j+s)}\right).\end{equation}
Setting $z=1$ recovers Selberg's zeta function. $\square$

\

Whenever we have a surface $X$ {and} refer to the flow-adapted IFS for $X$, we will assume the choice of potential function $V_s(u)=((f^{-1})'(u))^{-s}$ as above. Furthermore, we will use ${\L}^i_s$ as shorthand for ${\L}^i_{V_s}$, and ${\L}_s$ as shorthand for ${\L}_{V_s}$.

\section{Geometric Limits}
\label{geometriclimits}

In this section, we will establish control over the IFS functions $f_{ij}$ defined in the previous section.

\begin{definition}
\label{gooddef}

For $\eta>0$, a surface $X\in\M_{g,n}(L_1,\dots,L_n)$ is $\eta$-great if $\Phi(X)$ is trivalent, all corridors contain their intercostals, and if each corridor is cut into two half-corridors by its intercostal, the boundary arcs of the resulting half-corridors have length at least $\eta$.

\end{definition}

For the rest of this section, we will work with a fixed $\eta$-great surface $X\in \M_{g,n}(L_1,\dots,L_n)$, and assume the $L_h$ are $C$-bounded by some $L$.

\begin{remark}

This terminology comes from \cite{HT25}, although its usage here is very slightly adjusted: {being $\eta$-great in the sense of Definition \ref{gooddef} implies being $2\eta$-great in the sense of \cite{HT25}, while being $C_{g}\eta$-great in the sense of \cite{HT25} implies being $\eta$-great in the sense of Definition \ref{gooddef} for some constant $C_{g}=O_{g,n,C}(1)$.}

\end{remark}

\begin{figure}[h]

\centering
\includegraphics[scale=0.22]{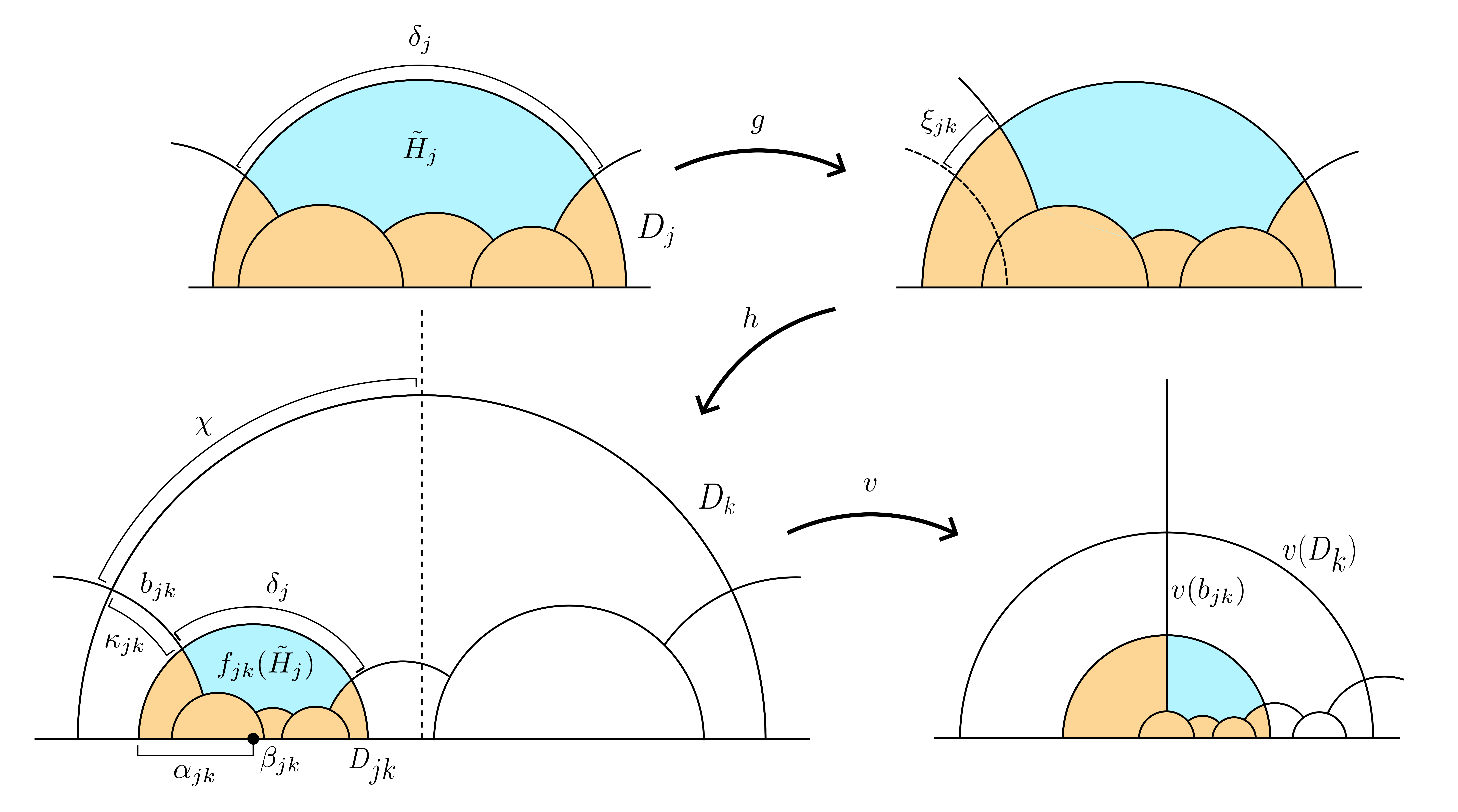}
\caption{An illustration of the maps, quantities, and objects used in the proof of Lemma \ref{mapcontrol}. Note that $f_{jk}=h\circ g$, where $g$ is a loxodromic isometry fixing the geodesic half-circle from $-1$ to $1$, while $h$ is an affine map. The map $v$ is an auxiliary isometry used in the proof of Lemma \ref{mapcontrol}. Note that $\alpha_{jk}$ is a Euclidean distance, while all other measured distances are hyperbolic distances.}

\label{kappafig}

\end{figure}

Assume $\Gamma=\Phi(X)$ has directed edges $e_1,\dots,e_{2{d}}$, and construct the flow-adapted IFS associated with $X$ as in subsection \ref{flowifs}. We want to obtain an explicit form for $f_{{jk}}:D_{{j}}\to D_{{jk}}$. To do so, we will decompose $f_{{jk}}$ into two functions, $f_{{jk}}=h\circ g$. The second function $h$ will be linear, of the form $h(u)=\alpha_{{jk}} u+\beta_{{jk}}$ for $\alpha_{{jk}},\beta_{{jk}}\in\reals$, and will be uniquely defined by the property $h(D_{{j}})=D_{{jk}}$. The second function will then be defined as the unique hyperbolic isometry fixing $D_{{j}}$ such that $f_{{jk}}=h\circ g$ as desired. It necessarily has the form
\begin{equation}g(u)=\frac{\cosh(\xi_{{jk}})u+\sinh(\xi_{{jk}})}{\sinh(\xi_{{jk}})u+\cosh(\xi_{{jk}})}\end{equation}
where $\xi_{{jk}}$ is its hyperbolic translation distance. {This setup is depicted in Figure \ref{kappafig}, which also depicts elements of the proof of Lemma \ref{mapcontrol}.}

\begin{figure}[h]

\centering
\includegraphics[scale=0.3]{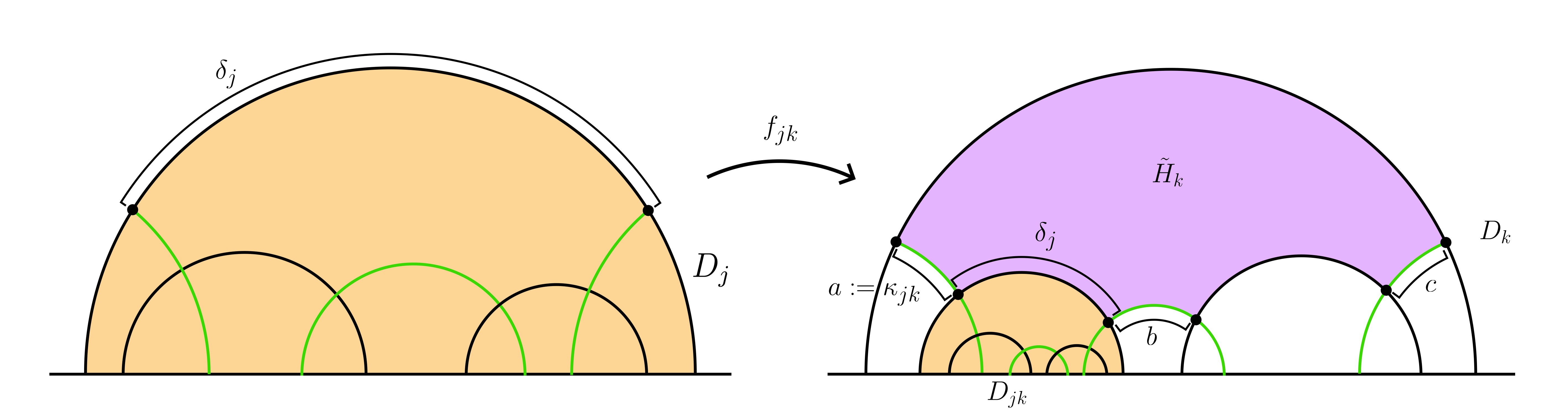}
\caption{The setup of Lemma \ref{intercostalcontrol}, showing the map $f_{jk}$ and the hexagon $\tilde{H}_k$ used in the proof. Note that this diagram is not to scale; in the situations where Lemma \ref{intercostalcontrol} applies, we expect $\kappa_{jk}$ to be very large, and $\delta_j$ to be very small.}

\label{deltafig}

\end{figure}

{Before moving on to the main estimate, we need a lemma. This lemma is illustrated in Figure \ref{deltafig}.}

\begin{lemma}
\label{intercostalcontrol}

Assume that $X\in\M_{g,n}(L_1,\dots,L_n)$ is $\eta$-great and that the $L_h$ are $C$-bounded by $L$. With the setup given above, for any $\delta_j$,
$$\delta_j=O(e^{-\eta}).$$

\end{lemma}

\textbf{Proof.} Choose an arbitrary $k$ so that $f_{jk}$ exists. Assume without loss of generality that $D_{{jk}}$ is the left {one} of the two disks contained in $D_{{k}}$. {Recall that $\delta_{{j}}=d_\H(\tilde{l}_{{j}},\tilde{r}_{{j}})$ and $\kappa_{{jk}}=d_\H(\partial D_{{k}},\partial D_{{jk}})$, where $\tilde{l}_j$ and $\tilde{r}_j$ are defined as in Figure \ref{sigmafig} and $D_k$ and $D_{jk}$ are defined as in Figure \ref{mapfig}.} {See Figure \ref{deltafig} for a complete illustration.}

To control $\delta_{{j}}$, label the boundary lengths of the half-wedges of the sector contained in $\tilde{H}_{{k}}$ as $a$, $b$, and $c$, so that $a=\kappa_{{jk}}$. This sector is a right-angled hexagon, so we may apply Lemma A.4 in \cite{GM09} to obtain
\begin{equation}\label{hexagon}\cosh(\delta_{{j}})=\frac{\cosh(a+b)\cosh(a+c)+\cosh(b+c)}{\sinh(a+b)\sinh(a+c)}.\end{equation}
Note the asymptotics $\cosh(x)=\frac{e^x}{2}+O(e^{-x})$, $\sinh(x)=\frac{e^x}{2}+O(e^{-x})$. Furthermore, essentially by definition of $\eta$-greatness, we have that $a,b,c\geq \eta$. Therefore, Equation \ref{hexagon} reduces to
$$\cosh(\delta_{{j}})=1+O(e^{-2\eta}).$$
We may now use $\cosh^2(x)-\sinh^2(x)=1$ to obtain $\sinh(\delta_{{j}})=O(e^{-\eta})$, and $\sinh(x)=x+O(x^3)$ to obtain $\delta_{{j}}=O(e^{-\eta})$. $\square$

{As mentioned previously, the following lemma is illustrated in Figure \ref{kappafig}.}

\begin{lemma}
\label{mapcontrol}

Assume that $X\in\M_{g,n}(L_1,\dots,L_n)$ is $\eta$-great and that the $L_h$ are $C$-bounded by $L$. With the setup given above, for any $f_{{jk}}$,
\begin{equation}
\label{mapeq}
f_{{jk}}(u)=\alpha_{{jk}}\left(\frac{\cosh(\xi_{{jk}})u+\sinh(\xi_{{jk}})}{\sinh(\xi_{{jk}})u+\cosh(\xi_{{jk}})}\right)+\beta_{{jk}}
\end{equation}
with
$$\Omega(e^{-L}){=} \alpha_{{jk}}{=}O(e^{-2\eta}),\ |\beta_{{jk}}|=O(e^{-\eta}),\ |\xi_{{jk}}|=O(e^{-\eta}).$$

\end{lemma}

\textbf{Proof.} Assume without loss of generality that $D_{{jk}}$ is the left {one} of the two disks contained in $D_{{k}}$.  {Recall as before that $\delta_{{j}}=d_\H(\tilde{l}_{{j}},\tilde{r}_{{j}})$ and $\kappa_{{jk}}=d_\H(\partial D_{{k}},\partial D_{{jk}})$, where $\tilde{l}_j$ and $\tilde{r}_j$ are defined as in Figure \ref{sigmafig} and $D_k$ and $D_{jk}$ are defined as in Figure \ref{mapfig}.} By construction, $\tilde{l}_j$ lies to the left of the imaginary axis, while $\tilde{r}_j$ lies to the right of this axis.

 Recall the decomposition $f_{jk}=h\circ g$ made above; {see Figure \ref{kappafig} for an illustration.} Let $\tilde{l}_{{j}}'=h^{-1}(f_{{jk}}(\tilde{l}_{{j}}))$ and $\tilde{r}_{{j}}'=h^{-1}(f_{{jk}}(\tilde{r}_{{j}}))$. Then $\xi=sd(\tilde{l}_{{j}},\tilde{l}_{{j}}')$, where $sd(\cdot,\cdot)$ is the signed hyperbolic distance from $\tilde{l}_{{j}}$ to $\tilde{l}_{{j}}'$ along the half-circle between $-1$ and $1$, with a clockwise orientation.

Drawing a vertical line through the Euclidean center of $D_{{jk}}$, it may be verified through visual inspection of {$f_{jk}(\tilde{H}_{{k}})$} {(recalling that $\tilde{H}_k$ is an interior region of the type depicted in Figure \ref{sigmafig})} that $f_{{jk}}(\tilde{l}_{{j}})$ must lie to the left of this line and $f_{{jk}}(\tilde{r}_{{j}})$ must lie to the right. This implies that $\tilde{l}_{{j}}'$ lies to the left of the imaginary axis, while $\tilde{r}_{{j}}'$ lies to the right, and in particular
$$|\xi_{{jk}}|=|sd(\tilde{l}_{{j}},\tilde{l}_{{j}}')|=d_\H(\tilde{l}_{{j}},\tilde{l}_{{j}}')\leq d_\H(\tilde{l}_{{j}},\tilde{r}_{{j}})+d_\H(\tilde{l}_{{j}}',\tilde{r}_{{j}}')=2\delta_{{j}}.$$

We may apply Lemma \ref{intercostalcontrol} to conclude $\delta_j =O(e^{-\eta})$.

We now handle $h$. By $h(D_{{j}})=D_{{jk}}$, $\beta_{{jk}}$ must be the (Euclidean) center of $D_{{jk}}$, while $\alpha_{{jk}}$ must be the radius of $D_{{jk}}$. To bound $\alpha_{{jk}}$, draw $\tilde{H}_{{k}}$ in $D_{{k}}$, and let $b_{{jk}}$ be the boundary arc realizing the distance from $\partial D_{{k}}$ to $\partial D_{{jk}}$. Earlier, we defined $\kappa_{{jk}}$ to be the length of this arc. We now apply a hyperbolic isometry $v$ {as in Figure \ref{kappafig}} to move $b_{{jk}}$ onto the imaginary axis while fixing $-1$ and $1$. This map sends $D_{{jk}}$ to the disk centered at $0$ with radius exactly $e^{-\kappa_{{jk}}}$. Furthermore, as $\kappa_{{jk}}$ is defined as the combined length of two incident half-wedge boundary arcs,
$$2\eta\leq \kappa_{{jk}}{=}O(L).$$

Explicitly,
$$v(u)= \frac{\cosh(\chi)u+\sinh(\chi)}{\sinh(\chi)u+\cosh(\chi)}$$
for some translation distance $\chi$ satisfying 
$$|\chi|=d_\H(\tilde{l}_{{k}},i)\leq d_\H(\tilde{l}_{{k}},\tilde{r}_{{k}})=\delta_{{k}}=O(e^{-\eta}).$$
{Here, the bound on $\delta_k$ follows as before from Lemma \ref{intercostalcontrol}.}
We therefore have that $\alpha_{{jk}}$, the Euclidean radius of $D_{{jk}}$, satisfies
$$\alpha_{{jk}}=\frac{|v^{-1}(e^{-\kappa_{{jk}}})-v^{-1}(-e^{-\kappa_{{jk}}})|}{2}.$$
But standard estimates and our bound on $|\chi|$ imply
$$v^{-1}(u)=\frac{(1+O(e^{-\eta}))u+O(e^{-\eta})}{O(e^{-\eta})u+{(}1+O(e^{-\eta}){)}}.$$
Applying this bound {plus our bounds on $\kappa_{jk}$ obtained above, and performing some straightforward simplifications, results in}
$$\Omega(e^{-L}){=} \alpha_{ij}{=} O(e^{-2\eta}).$$

To bound $\beta_{{jk}}$, a straightforward argument shows that $|\beta_{{jk}}|$ is bounded above by the Euclidean distance from $i$ to $\tilde{l}_{{k}}$. But we may argue as above to show $d_\H(i,\tilde{l}_{{k}})=O(e^{-\eta})$. For sufficiently large $\eta$, $d_\complex(i,\tilde{l}_{{k}}){=O(} d_\H(i,\tilde{l}_{{k}}){)}$ (as the hyperbolic metric and the Euclidean metric approximately agree near $i$), so $|\beta_{{jk}}|=O(e^{-\eta})$. $\square$

\

Recall that $D=\bigcup_{{j}=1}^{2{d}}D_{{j}}$, and consider some arbitrary $a\in A_\infty(D)$. Each $f_{{jk}}$ is notionally only defined on $D_{{j}}$; however, as a linear fractional transformation, it extends uniquely to a meromorphic function on all of $\complex$. As a consequence, $a\circ f_{{jk}}$ can be extended uniquely to a function on $f_{{jk}}^{-1}(D_{{k}})$. For a later step, it will be crucial to understand the maximal radius around $0$ on which this extension is defined.

\begin{lemma}
\label{inversecontrol}

Let $X\in\M_{g,n}(L_1,\dots,L_n)$ be $\eta$-great {with the $L_h$ $C$-controlled by $L$. {Let $\epsilon>0$,} and assume $\eta$ is sufficiently large in terms of $g$, $n$, $C$, {and $\epsilon$.}}  With the previously established notation, for all functions $f_{{jk}}$ defined as part of the flow-adapted IFS corresponding to $X$, $f_{{jk}}^{-1}(D_{{k}})$ contains the ball of radius {$e^{(1-\epsilon)\eta}$} centered at $0$ in its interior.
	
\end{lemma}

\textbf{Proof.} It suffices to prove that
$${e^{(1-\epsilon)\eta}}\leq \min(|f_{{jk}}^{-1}(-1)|,|f_{{jk}}^{-1}(1)|).$$
Using equation \ref{mapeq}, we may compute that
\begin{equation}f_{{jk}}^{-1}(u)=\frac{\cosh(-\xi_{{jk}})(\frac{1}{\alpha_{{jk}}}(u-\beta_{{jk}}))+\sinh(-\xi_{{jk}})}{\sinh(-\xi_{{jk}})(\frac{1}{\alpha_{{jk}}}(u-\beta_{{jk}}))+\cosh(-\xi_{{jk}})}.\end{equation}
By lemma \ref{mapcontrol},  we have the bounds
$$0<\alpha_{{jk}}<e^{-2\eta},\ |\beta_{{jk}}|{=}O(e^{-\eta}),\ |\xi_{{jk}}|{=}O(e^{-\eta}).$$
We will consider $u=1$; the case of $u=-1$ will be analogous. Applying asymptotics,
$$f_{{jk}}^{-1}(1)=\frac{(1+O(e^{-2\eta}))(\frac{1}{\alpha_{{jk}}}(1+O(e^{-\eta})))+O(e^{-\eta})}{O(e^{-\eta})(\frac{1}{\alpha_{{jk}}}(1+O(e^{-\eta})))+1+O(e^{-2\eta})}$$
\begin{equation}=\frac{\frac{1}{\alpha_{{jk}}}(1+O(e^{-\eta}))}{\frac{1}{\alpha_{{jk}}}O(e^{-\eta})+O(1)}.\end{equation}
Multiplying through by $\alpha_{{jk}}$ and using ${0<\alpha_{jk}<e^{-2\eta}}$ yields
$${f^{-1}_{jk}(1)=}\frac{1+O(e^{-\eta})}{O(e^{-\eta})}=\Omega(e^\eta)$$
and so we conclude that
$$f_{{jk}}^{-1}(1)=\Omega(e^{\eta}).$$
{For sufficiently large $\eta$ in terms of $g$, $n$, $C$, {and $\epsilon$,} this bound implies that ${e^{(1-\epsilon)\eta}}\leq f_{jk}^{-1}(1)$.}  $\square$

\

Lastly, we need some control on derivatives of the $f_{{jk}}$.

\begin{lemma}
\label{derivcontrol}

Let $X\in \M_{g,n}(L_1,\dots,L_n)$ be $\eta$-great {with the $L_h$ $C$-controlled by $L$. {Also let $\epsilon>0$ and} assume $\eta$ is sufficiently large in terms of $g$, $n$, $C$, { and $\epsilon$}}. For all functions $f_{{jk}}$ defined as part of the flow-adapted IFS corresponding to $X$, and for all $|u|\leq {e^{(1-\epsilon)\eta}}$,
$$\Omega(e^{-L}){=}|f'_{{jk}}(u)|{=}O(e^{-\eta}).$$
	
\end{lemma}

\textbf{Proof.} By Equation \ref{mapeq}, we have
\begin{equation}f'_{{jk}}(u)=\frac{\alpha_{{jk}}}{(\sinh(\xi_{{jk}})u+\cosh(\xi_{{jk}}))^2}.\end{equation}
Applying Lemma \ref{mapcontrol}, this yields
\begin{equation}f'_{{jk}}(u)=\frac{\alpha_{{jk}}}{(O(e^{-\eta})u+1+O(e^{-\eta}))^2}.\end{equation}
For $|u|\leq {e^{(1-\epsilon)\eta}}$ and for $\eta$ sufficiently large, the denominator is $({O(e^{-\epsilon\eta})}+1+O(e^{-\eta}))^2$, and the desired bounds follow from the bounds on $\alpha_{{jk}}$ obtained in Lemma \ref{mapcontrol}. $\square$

\

For the functions $f_{{jk}}$ as defined in Lemma \ref{derivcontrol}, we may therefore guarantee that, say, $|f'_{{jk}}(u)|\leq\frac{1}{2}$ for all $u$ in the domain of $f_{{jk}}$ as long as $\eta$ is sufficiently large. The same bound holds for any $|f'_w(u)|$ by the chain rule, so we obtain the following:

\begin{corollary}

Let $X\in\M_{g,n}(L_1,\dots,L_n)$ be $\eta$-great {with the $L_h$ $C$-controlled by $L$, and assume $\eta$ is sufficiently large in terms of $g$, $n$, and $C$}. Then the flow-adapted IFS corresponding to $X$ is eventually contracting.
	
\end{corollary}

\section{The Nuclear Representation and Cycle Expansion}

\subsection{Controlling the Nuclear Representation}

Recall that transfer operators of eventually contracting holomorphic IFSes are known to be nuclear operators. In a sense, this property is the entire motivation for rewriting the Selberg zeta function as a dynamical zeta function, as it provides access to the following powerful result:

\begin{lemma}[Grothendieck \cite{AG56}] 
\label{nuclear}

If $B$ is a Banach space and $L:B\to B$ is a nuclear operator with the nuclear representation
\begin{equation}{\L}h=\sum_{{m}=0}^\infty \lambda_{{m}}a_{{m}}(h)v_{{m}}.\end{equation}
Then the Fredholm determinant $\det(1-z {\L})$ has a power series expansion
$$1+\sum_{{m}=0}^\infty z^{{m}}{c}_{{m}}$$
defined for all $z\in\complex$ and satisfying
\begin{equation}
\label{cycleeq}
{c}_{{m}}=(-1)^{{m}}\sum_{{q}_1<\cdots<{q}_k}\lambda_{{q}_1}\cdots\lambda_{{q}_k}\det\left((a_{{q}_{{m_1}}}(v_{{q}_{{m_2}}}))_{{m_1},{m_2}=1}^{{m}}\right)
\end{equation}
where $(a_{{q}_{{m_1}}}(v_{{q}_{{m_2}}}))_{{m_1},{m_2}=1}^{{m}}$ is the ${m}\times {m}$ matrix with entries $\alpha_{{q}_{{m_1}}}(v_{{q}_{{m_2}}})$.

\end{lemma}

The power series expansion $\det(1-z{\L})=1+\sum_{{m}=0}^\infty z^{{m}}{c}_{{m}}$ given above is known as the \emph{cycle expansion} of the zeta function.

We will now discuss how a nuclear decomposition may be found for our setup. Note that since our operator is technically a family of operators parameterized by $s$, the terms in the nuclear decomposition will likewise be parameterized by $s$.

 The starting point is the operators ${\L}_s^{{j}}:A_\infty(D)\to A_\infty(D_{{j}})$; recall that ${\L}_s^{{j}}$ is shorthand for ${\L}^{{j}}_{V_s}$ with $V_s$ defined as in Lemma \ref{selbergequality}. We may extend the range of each ${\L}_s^{{j}}$ to $A_\infty(D)$ by setting it equal to $0$ on all other $A_\infty(D_{{k}})$, and then define ${\L}_s$ as the sum of these extended operators. In particular, we can calculate the nuclear representations for each ${\L}^{{j}}_s$ separately, and add them to obtain the nuclear representation for ${\L}_s$.

\begin{lemma}
\label{nuclearcontrol}

Let $X\in\M_{g,n}(L_1,\dots,L_n)$ be $\eta$-great with the $L_{{h}}$ $C$-bounded by $L$, and let the ${\L}_s^{{j}}$ be the components of the transfer operator for the flow-adapted IFS corresponding to $X$. {Also let $\epsilon>0$.} {Assume that $\eta$ is sufficiently large in terms of $g$, $n$, $C$, and {$\epsilon$}.}  Then each ${\L}^{{j}}_s$ has the nuclear representation
\begin{equation}{\L}^{{j}}_s=\sum_{{m}=0}^\infty \lambda^{{j}}_{{m},s}a_{{m},s}^{{j}}v^{{j}}_{{m},s}\end{equation}
with terms satisfying the formulas
	
$$v^{{j}}_{{m},s}(u)=u^{{m}},$$
$$\tilde{a}^{{j}}_{{m},s}(h)=\frac{1}{2\pi i}\int_{C_\rho}\frac{{\L}_s^{{j}}h(w)}{w^{{m}+1}}dw,$$
$$a_{{m},s}^{{j}}(h)=\tilde{a}^{{j}}_{{m},s}(h)/\|\tilde{a}^{{j}}_{{m},s}\|,$$
\begin{equation}|\lambda_{{m},s}^{{j}}|={O(e^{-{m}{(1-\epsilon)}\eta}e^{L|\re(s)|})}.\end{equation}
Here, $C_\rho$ is the circle around the origin with radius $\rho={e^{(1-\epsilon)\eta}}$.

\end{lemma}

\textbf{Proof.} Fix ${\L}^{{j}}_s$, and recall the definition
\begin{equation}{\L}^{{j}}_sh(u)=\sum_{{k}:{j}\to {k}}V_s(f_{{jk}}(u))h(f_{{jk}}(u)).\end{equation}
As was argued earlier, $f_{{jk}}$ maps a disk centered at $0$ with radius $\rho={e^{(1-\epsilon)\eta}}$ into $D_{{k}}$, and is holomorphic on a slightly larger disk (in fact, a much larger disk, with radius $\Omega(e^\eta)$). Let $C_\rho$ be the circle around $0$ with radius $\rho$. By the Cauchy integral formula, for any $u\in D_{{j}}$,
\begin{equation}{\L}^{{j}}_sh(u)=\frac{1}{2\pi i}\int_{C_\rho}\frac{{\L}^{{j}}_sh(w)}{w-u}dw.\end{equation}
We rewrite this expression as
\begin{equation}\frac{1}{2\pi i}\int_{C_\rho}\frac{{\L}_s^{{j}}h(w)}{w}\frac{1}{1-\frac{u}{w}}=\sum_{{m}=0}^\infty\frac{1}{2\pi i}\int_{C_\rho}\frac{{\L}^{{j}}_sh(w)}{w}\left(\frac{u}{w}\right)^ndw.\end{equation}
We may therefore define
$$\tilde{a}_{{m},s}^{{j}}(h)=\frac{1}{2\pi i}\int_{C_\rho}\frac{{\L}^{{j}}_sh(w)}{w^{{m}+1}}dw,$$
\begin{equation}\tilde{v}_{{m},s}^{{j}}(u)=u^{{m}},\end{equation}
leading to
$$a_{{m},s}^{{j}}=\frac{\tilde{a}_{{m},s}^{{j}}}{\|\tilde{a}^{{j}}_{{m},s}\|},\ v^{{j}}_{{m},s}(u)=\frac{\tilde{v}_{{m},s}^{{j}}}{\|\tilde{v}_{{m},s}^{{j}}\|},$$
\begin{equation}\lambda_{{m},s}^{{j}}=\|\tilde{a}^{{j}}_{{m},s}\|\|\tilde{v}^{{j}}_{{m},s}\|.\end{equation}
These terms give a valid nuclear decomposition for ${\L}^{{j}}_s$ (we will check in a moment that the $\lambda^{{j}}_{{m},s}$ satisfy the bounded series condition). 

We would like to obtain estimates on the $\lambda^{{j}}_{{m},s}$. Since the norm on $\tilde{v}^{{j}}_{{m},s}$ is the supremum norm on $D_{{j}}$, $\|\tilde{v}^{{j}}_{{m},s}\|=\|u^{{m}}\|$ is clearly $1$. For $\|\tilde{a}_{{m},s}^{{j}}\|$, we may take arbitrary $h\in A_\infty(D)$ satisfying $\|h\|\leq 1$ and apply a supremum bound:
\begin{equation}|\tilde{a}_{{m},s}^{{j}}(h)|=\left|\frac{1}{2\pi i}\int_{C_\rho}\frac{{\L}_s^{{j}} h(w)}{w^{{m}+1}}dw\right|\leq \rho^{-{m}}\|{\L}^{{j}}_sh\|_{L^\infty(C_\rho)}.\end{equation}
Recall that
\begin{equation}{\L}^{{j}}_sh(u)=\sum_{{k}:{j}\to {k}}V_s(f_{{jk}}(u))h(f_{{jk}}(u)).\end{equation}
By trivalency, the sum consists of two terms, so it suffices to bound each term individually. Furthermore, $|h(f_{{jk}}(u))|\leq 1$ by assumption. On the other hand,
$$V_s(f_{{jk}}(u))=(f^{-1}(f_{{jk}}(u))^{-s}=f_{{jk}}'(u)^s.$$
For $u\in C_\rho$, we have by Lemma \ref{derivcontrol} that $|f'_{{jk}}(u)|{^s}=O(e^{L|\re(s)|})$. Therefore
$$\|{\L}^{{j}}_sh\|_{L^\infty(C_\rho)}=O(e^{L{|\re(s)|}}).$$

As a result, {we may plug in $\rho= {e^{(1-\epsilon)\eta}}$ and obtain}
$$|\lambda^{{j}}_{{m},s}|={O(e^{-{m}{(1-\epsilon)\eta}}e^{L{|\re(s)|}})}.\text{ $\square$}$$

\

We therefore have the following nuclear decomposition of ${\L}_s$:
\begin{equation}{\L}_sh=\sum_{{j}=1}^{2{d}}\sum_{{m}=0}^\infty\lambda^{{j}}_{{m},s}a^{{j}}_{{m},s}(h)v^{{j}}_{{m},s}\end{equation}
and the exponential decay of the $\lambda^{{j}}_{{m},s}$ in ${m}$ obtained as part of Lemma \ref{nuclearcontrol} guarantees that $\sum_{{j}=1}^{2{d}}\sum_{{m}=0}^\infty |\lambda^{{j}}_{{m},s}|<\infty$ for any fixed $s$.

\subsection{Bounding the Tail}

\begin{lemma}
\label{truncation}

Let $X\in\M_{g,n}(L_1,\dots,L_n)$ be $\eta$-great, assume the $L_{{h}}$ are $C$-bounded by $L$, and let ${\L}_s$ be the transfer operator for the flow-adapted IFS corresponding to $X$. Fix some constant $A$, and complex numbers $s,z$ with $|z|\leq A$, $|s|\leq \frac{A}{L}$. {Let $\epsilon>0$, and} assume that $\eta$ is sufficiently large in terms of $g$, $n$, $C$, $A$, and {$\epsilon$}. Then the cycle expansion
\begin{equation}\det(1-z{\L}_s)=1+\sum_{{m}=1}^\infty z^{{m}}{c}_{{m}}(s)\end{equation}
satisfies the truncation bound
\begin{equation}\left|\sum_{{m}=2{d}+1}^\infty z^{{m}}{c}_{{m}}(s)\right|{=O_{g,n,C,A,{\epsilon}}(e^{-{(1-\epsilon)\eta}})}.\end{equation}
	
\end{lemma}

\textbf{Proof.} {Note that our eventual goal is to show we obtain a convergent geometric series, so we need to be somewhat careful with implicit constants.} 

We begin from Formula \ref{cycleeq}, and wish to use Lemma \ref{nuclearcontrol} to understand the cycle expansion coefficients. To start, fix ${m}$ and consider ${c}_{{m}}(s)$. This term is defined via a sum indexed by $2{d}$-tuples of increasing sequences of nonnegative integers
$$({p}_{1,1}<{p}_{2,1}<\cdots<{p}_{r_1,1},{p}_{1,2}<{p}_{2,2}<\cdots <{p}_{r_2,2},\dots,{p}_{1,2k}<\cdots<{p}_{r_{2{d}},2{d}})$$
satisfying $r_1+\cdots+r_{2{d}}={m}$. Take such a tuple. Its contribution consists of a product of the determinant of a ${m}\times {m}$ matrix (with all entries bounded by $1$) and the expression
$$\prod_{{j}=1}^{2{d}}\prod_{{q}=0}^{r_{{j}}}\lambda^{{j}}_{{p}_{{q},{j}},{s}}.$$
Notice that we allow some $r_{{j}}$ to be $0$ in this notation. 

By the Hadamard bound, the determinant component is bounded in absolute value by ${m}^{{m}/2}$. Furthermore, by Lemma \ref{nuclearcontrol}, {there is some constant $C_1$ dependent only on $g$, $n$, $C$, $A$, {and $\epsilon$} so that}
$$\prod_{{j}=1}^{2{d}}\prod_{{q}=0}^{r_{{j}}}|\lambda^{{j}}_{{p}_{{q},{j}},s}|\leq\prod_{{j}=1}^{2{d}}\prod_{{q}=0}^{r_{{j}}}{C_1}e^{-{p}_{{q},{j}}{(1-\epsilon)\eta}+L|\re(s)|}$$
\begin{equation}=C_1^me^{{m}L|\re(s)|}\prod_{{j}=1}^{2{d}}\prod_{{q}=0}^{r_{{j}}}e^{-{p}_{{q},{j}}{(1-\epsilon)\eta}}.\end{equation}
Therefore, we have the following bound for ${m}\geq 1$:
\begin{equation}|{c}_{{m}}(s)|{\leq}{C_1^m}  {m}^{{m}/2}e^{{m}L|\re(s)|}\sum_{*}\prod_{{j}=1}^{2{d}}\prod_{{q}=0}^{r_{{j}}}e^{-{p}_{{q},{j}}{(1-\epsilon)\eta}}\end{equation}
where the sum is over all tuples
$$({p}_{1,1}<\cdots<{p}_{r_1,1},\dots,{p}_{1,2k}<\cdots<{p}_{r_{2{d}},2{d}})$$
such that $r_1+\cdots+r_{2{d}}={m}$. To handle the sum, we first split it into two, fixing $r_1,\dots,r_{2{d}}$ in the first sum:
\begin{equation}{\sum_{*}\prod_{{j}=1}^{2{d}}\prod_{{q}=0}^{r_{{j}}}e^{-{p}_{{q},{j}}{(1-\epsilon)\eta}}=}\sum_{\substack{r_1,\dots,r_{2{d}}:\\r_1+\cdots+r_{2{d}}={m}}}\sum_{*'}\prod_{{j}=1}^{2{d}}\prod_{{q}=0}^{r_{{j}}}e^{-{p}_{{q},{j}}{(1-\epsilon)\eta}}.\end{equation}
This expression becomes
\begin{equation}\sum_{\substack{r_1,\dots,r_{2{d}}:\\r_1+\cdots+r_{2{d}}={m}}}\prod_{{j}=1}^{2{d}}\sum_{{p}_{1,{j}}<\cdots<{p}_{r_{{j}},{j}}}e^{-({p}_{1,{j}}+\cdots+{p}_{r_{{j}},{j}}){(1-\epsilon)\eta}}.\end{equation}
Let $\theta_{r}({t})$ count the number of strictly increasing sequences ${p}_{1}<\cdots<{p}_{r}$ satisfying ${p}_1+\cdots+{p}_r={t}$, so that
\begin{equation}\sum_{{p}_{1,{j}}<\cdots<{p}_{r_{{j}},{j}}}e^{-({p}_{1,{j}}+\cdots+{p}_{r_{{j}},{j}})\eta/2}=\sum_{{t}=0}^\infty\theta_{r_{{j}}}({t})e^{-{t}{(1-\epsilon)\eta}}.\end{equation}

The quantity $\theta_r({t})$ is bounded by the number of $k$-partite partitions ${p}_1+\cdots+{p}_r={t}$ (without the strictly increasing condition), which is in turn bounded above by  ${t}^r$. Furthermore, {the strictly increasing condition implies that} $\theta_{{r}}({t})=0$ for ${t}< {\frac{r^2+r}{2}:=T(r)}$. So
\begin{equation}|{c}_{{m}}(s)|\leq {C_1^m} {m}^{{m}/2}e^{{m}L|\re(s)|}\sum_{{r}_1+\cdots+{r}_{2{d}}={m}}\prod_{{j}=1}^{2{d}}\sum_{{t}={T(r_j)}}^\infty {t}^{r_{{j}}}e^{-{t}{(1-\epsilon)\eta}}.\end{equation}

{To handle the general expression
$$\sum_{t=T(r)}^\infty t^r e^{-t{(1-\epsilon)\eta}}$$
we first rewrite it as
$$\sum_{t=0}^\infty (t+T(r))^re^{-(t+T(r)){(1-\epsilon)\eta}}$$
\begin{equation}=e^{-T(r)\eta/2}\sum_{r'=0}^r\binom{r}{r'}T(r)^{r-r'}\sum_{t=0}^\infty t^{r'}e^{-t {(1-\epsilon)\eta}}\end{equation}
We can uniformly bound $\binom{r}{r'}T(r)^{r-r'}$ by $r^{4r}$, while $\sum_{t=0}^\infty t^{r'}e^{-t{(1-\epsilon)\eta}}$ is by standard Gamma function estimates at most $2\frac{r!}{({(1-\epsilon)\eta})^{r+1}}\leq 2\frac{r^r}{({(1-\epsilon)\eta})^{r+1}}$. For our purposes, the bound $\sum_{t=0}^\infty t^{r'}e^{-t{(1-\epsilon)\eta}}\leq 2r^r$ is sufficient.}

{Returning to our previous bound on $|c_m(s)|$, we obtain
\begin{equation}|{c}_{{m}}(s)|\leq {C_2^m} {m}^{{m}/2}e^{{m}L|\re(s)|}\sum_{{r}_1+\cdots+{r}_{2{d}}={m}}\prod_{{j}=1}^{2{d}}r_j^{4r_j}e^{-T(r_j){(1-\epsilon)\eta}}\end{equation}
where $C_2$ is a different constant depending only on $g$, $n$, $C$, $A$, {and $\epsilon$}.}

{Applying the bound $r_j\leq m$ leads to a bound of
$$C_2^m m^{m/2}e^{mL|\re(s)|}\sum_{r_1+...+r_{2d}=m}\prod_{j=1}^{2d}m^{4r_j}e^{-T(r_j){(1-\epsilon)\eta}}$$
\begin{equation}\label{midbound}\leq C_2^m m^{5m}e^{mL|\re(s)|}\sum_{r_1+...+r_{2d}=m}\prod_{j=1}^{2d}e^{-T(r_j){(1-\epsilon)\eta}}\end{equation}
We now note that for any nonnegative integers $r_1+...+r_{2d}=m$ we have that some $r_j\geq \frac{m}{2d}$, and so
$$T(r_1)+...+T(r_{2d})\geq T\left(\frac{m}{2d}\right)\geq \frac{m^2}{8d^2}.$$
On the other hand, $T(r)\geq r$ for all nonnegative $r$, so
$$T(r_1)+...+T(r_{2d})\geq m$$
We will merge these inequalities by defining a piecewise function
$$\mathcal{T}(m)=\begin{cases} m & m\leq M \\ \frac{m^2}{8d^2} & m> M\end{cases}$$
dependent on some quantity $M$ that will be chosen later in the proof. For all $r_1+...+r_{2d}=m$, we have $T(r_1)+...+T(r_{2d})\geq \mathcal{T}(m)$ by construction, leading to an upper bound for Equation \ref{midbound} of
\begin{equation}\leq C_2^m m^{5m}e^{mL|\re(s)|}e^{-\mathcal{T}(m){(1-\epsilon)\eta}}\sum_{r_1+...+r_{2d}=m}1.\end{equation}
A standard combinatorial fact states that the number of ordered nonnegative tuples $(r_1,...,r_{2d})$ such that $r_1+...+r_{2d}=m$ is $\binom{m+2d-1}{2d-1}$, and by standard bounds
$$\binom{m+2d-1}{2d-1}\leq C_3 m^{2d}$$
where $C_3:=C_3(d)$ is a constant depending only on $d$. We thus obtain that
\begin{equation}|c_m(s)|\leq C_4^m m^{5m+2d}e^{mL|\re(s)|}e^{-\mathcal{T}(m){(1-\epsilon)\eta}}\end{equation}
for some constant $C_4$ depending only on $g$, $n$, $C$, $A$, {and $\epsilon$}.
Recall that our goal is to bound
$$\left|\sum_{m=2d+1}^\infty z^m c_m(s)\right|.$$
Applying $|z|\leq A$ and $L|\re(s)|\leq A$ by the assumptions on $z$ and $s$, we have
$$\left|\sum_{m=2d+1}^\infty z^m c_m(s)\right|\leq \sum_{m=2d+1}^\infty |z|^m|c_m(s)|$$
$$\leq \sum_{m=2d+1}^\infty A^mC_4^mm^{5m+2d}e^{Am}e^{-\mathcal{T}(m){(1-\epsilon)\eta}}$$
\begin{equation}\leq\sum_{m=2d+1}^\infty\left(\frac{(C_4Ae^A) m^{6}}{e^{\mathcal{T}(m){(1-\epsilon)\eta}m}}\right)^m.\end{equation}
Recall the $M$ used to define $\mathcal{T}(m)$, and consider
\begin{equation}\sum_{m=M}^\infty\left(\frac{(C_4Ae^A)m^6}{e^{\mathcal{T}(m){(1-\epsilon)\eta}m}}\right)^m\leq \sum_{m=M}^\infty\left(\frac{(C_4Ae^A)m^6}{e^{m{(1-\epsilon)\eta}/4d^2}}\right)^m.\end{equation}
Standard asymptotics show that for $\eta$, say, at least $10$, the expression
$$\frac{(C_4Ae^A)m^6}{e^{m{(1-\epsilon)\eta}/4d^2}}$$
is strictly decreasing for all $m$ sufficiently large in terms of $g$, $n$, $C$, $A$, {and $\epsilon$} (and hence sufficiently large in terms of $d$). As a result, there is some $M_0$ dependent only on $g$, $n$, $C$, $A$, {and $\epsilon$} so that for $M\geq M_0$, we have by a standard geometric series bound that
$$\sum_{m=M}^\infty\left(\frac{(C_4Ae^A)m^6}{e^{m{(1-\epsilon)\eta}/4d^2}}\right)^m\leq\sum_{m=M}^\infty\left(\frac{(C_4Ae^A)M^6}{e^{M{(1-\epsilon)\eta}/4d^2}}\right)^m$$
\begin{equation}\label{tailgeom}\leq\frac{(C_4Ae^A)M^{6M}}{e^{M^2{(1-\epsilon)\eta}/4d^2}}*\frac{\frac{(C_4Ae^A)M^6}{e^{M{(1-\epsilon)\eta}/4d^2}}}{1-\frac{(C_4Ae^A)M^6}{e^{M{(1-\epsilon)\eta}/4d^2}}}\end{equation}
We now can find $M\geq M_0$ (this will be our chosen value of $M$ for the rest of the proof) dependent only on $g$, $n$, $C$, $A$, and $\epsilon$, and $\eta_0$ dependent only on $g$, $n$, $C$, $A$, and $\epsilon$, so that for all $\eta\geq \eta_0$,
\begin{enumerate}
\item $\frac{(C_4Ae^A)M^6}{e^{M{(1-\epsilon)\eta}/4d^2}}\leq \frac{1}{2}$.
\item $\frac{M^2}{{4d^2}}\geq {\frac{1}{1-\epsilon}}$.\end{enumerate}
With these inequalities, for all $\eta\geq \eta_0$, Equation \ref{tailgeom} is $O_{g,n,C,A,{\epsilon}}(e^{-\eta})$
We now turn to
\begin{equation}\sum_{m=2d+1}^M \left(\frac{(C_4Ae^A)m^6}{e^{\mathcal{T}(m){(1-\epsilon)\eta}m}}\right)^m.\end{equation}
By definition of $\mathcal{T}(m)$ and $M$, this expression is bounded above by
\begin{equation}\sum_{m=2d+1}^M\left(\frac{(C_4Ae^A)m^6}{e^{{(1-\epsilon)\eta}}}\right)^m.\end{equation}
In turn, we can use $1\leq m\leq M$ to apply a simple bound of
\begin{equation}\leq M\frac{(C_rAe^A)^MM^{6M}}{e^{{(1-\epsilon)\eta}}}=O_{g,n,C,A,{\epsilon}}(e^{-{(1-\epsilon)\eta}})\end{equation}
Combining our bounds, we conclude that
\begin{equation}\left|\sum_{m=2d+1}^\infty z^m c_m(s)\right|=O_{g,n,C,A,{\epsilon}}(e^{-{(1-\epsilon)\eta}}).\ \square\end{equation}

\section{Convergence of Non-Decaying Terms}

Having controlled $\sum_{{m}=2{d}+1}^\infty z^{{m}} {c}_{{m}}(s)$, we must now handle the rest of the expression. We start with an alternate formula for the cycle expansion terms.

\begin{lemma}[Jenkinson-Pollicott, Proposition 8 \cite{JP02}]
\label{altformula}

Fix $X\in\M_{g,n}(L_1,\dots,L_n)$, and let ${\L}_s$ be the transfer operator for the flow-adapted IFS corresponding to $X$. Then the terms of the cycle expansion
\begin{equation}\det(1-z{\L}_s)=1+\sum_{{m}=1}^\infty z^{{m}} {c}_{{m}}(s)\end{equation}
satisfy the formula
\begin{equation}{c}_{{m}}(s)=\sum_{{r}=1}^{{m}}\sum_{({p}_1,\dots,{p}_{{r}})\in P({m},{r})}\frac{(-1)^{{r}}}{{r}!}\prod_{{q}=1}^{{r}}\frac{1}{{p}_{{q}}}\sum_{w\in W^C_{{p}_{{q}}}}\frac{e^{-sl(\gamma_w)}}{1-e^{-l(\gamma_w)}}\end{equation}
where $P({m},{t})$ is the set of all ${t}$-fold partitions of ${m}$.

\end{lemma}

An analogous formula holds for graphs:

\begin{lemma}
\label{altformulagraphs}

Let $\Gamma$ be a finite {(undirected)} metric graph with ${d}$ edges, and construct $J_s$ as in Subsection \ref{zetafunctions}. For ${p}\geq 1$, let $ W^C_{{p}}$ be the set of all closed words {with length $p$} in the directed-edge symbolic coding. Then we have the formula
\begin{equation}\det(1-zJ_s)=1+\sum_{{m}=1}^{2{d}}z^{{m}}\tilde{{c}}_{{m}}(s)\end{equation}
with
\begin{equation}\tilde{{c}}_{{m}}(s)=\sum_{{r}=1}^{{m}}\sum_{({p}_1,\dots,{p}_{{r}})\in P({m},{r})}\frac{(-1)^{{r}}}{{r}!}\prod_{{q}=1}^{{r}}\frac{1}{{p}_{{q}}}\sum_{w\in W^C_{{p}_{{q}}}}e^{-sl(\tilde{\gamma}_w)}.\end{equation}
	
\end{lemma}

\textbf{Proof.} {We will only give a sketch, as the proof is very similar to that of Lemma \ref{altformula}.} As a matrix, $J_s$ has dimension $2{d}\times 2{d}$. So the cycle expansion is a polynomial in $z$ with degree at most $2{d}$. To obtain the exact formula, start from the identity
\begin{equation}\det(1-zJ_s)=\exp\left(-\sum_{{m}=1}^\infty\frac{z^{{m}}}{{m}}\tr(J_s^{{m}})\right)\end{equation}
and observe that $\tr(J_s^{{m}})=\sum_{w\in W^C_{{m}}}e^{-sl(\tilde{\gamma}_w)}$. Then take a power series expansion and group terms. $\square$

\

We aim to compare these two types of terms.

\begin{lemma}
\label{earlyconvergence}

Let $X\in\M_{g,n}(L_1,\dots,L_n)$ be $\eta$-great with the $L_h$ $C$-bounded by $L$, and set $\Gamma=\Phi(X)$. Let ${\L}_s$ be the transfer operator derived from the flow-adapted IFS corresponding to $X$, and let $J_s$ be the transfer operator corresponding to $\Gamma$. Define the cycle expansion terms ${c}_{{m}}(s)$ and $\tilde{{c}}_{{m}}(s)$ as in Lemmas \ref{altformula} and \ref{altformulagraphs}. Fix a particular ${m}\leq 2{d}$, and $z,s\in\complex$ satisfying $|z|\leq A$, $|s|\leq \frac{A}{L}$ for some chosen constant $A$. Then
\begin{equation}|{c}_{{m}}(s)-\tilde{{c}}_{{m}}(s)|= {O_{g,n,C,A}(e^{-\eta})}.\end{equation}
	
\end{lemma}

\textbf{Proof.} From Lemma 4.20 in \cite{HT25} (and after accounting for our slightly different definition of $\eta$-greatness), we know that for all $\gamma_w$ on $X$ with word length ${m}$,
$$l(\tilde{\gamma}_w)\leq l(\gamma_w)\leq l(\tilde{\gamma}_w)+O({m}e^{-\eta}).$$
Now consider
$$\frac{e^{-sl(\gamma_w)}}{1-e^{-l(\gamma_w)}}-e^{-sl(\tilde{\gamma}_w)}$$
\begin{equation}=\frac{e^{-sl(\gamma_w)}-(1-e^{-l(\gamma_w)})e^{-sl(\tilde{\gamma}_w)}}{1-e^{-l(\gamma_w)}}\end{equation}
for some word $w$ of length ${m}$. Since both $l(\gamma_w)$ and $l(\tilde{\gamma}_w)$ are bounded below by $\eta$, for $\eta$ absolutely large, the denominator is between, say, $1/2$ and $2$. Meanwhile, the numerator expands to
\begin{equation}e^{-sl(\gamma_w)}-e^{-sl(\tilde{\gamma}_w)}+e^{-sl(\tilde{\gamma}_w)-l(\gamma_w)}.\end{equation}
This expression factors to
\begin{equation}e^{-sl(\gamma_w)}(1-e^{-s(l(\tilde{\gamma}_w)-l(\gamma_w))}(1+e^{-l(\gamma_w)}))\end{equation}
\begin{equation}=e^{-sl(\gamma_w)}(1-e^{-sO({m}e^{-\eta})}(1+O(e^{-\eta}))).\end{equation}
Since $l(\gamma_w)\leq O(L{m})$ and $|s|\leq \frac{A}{L}$, we may apply a standard Taylor estimate and obtain a bound on the preceding expression of
$$O(e^{A{m}})(1-(1+O(A{m}e^{-\eta}/L))(1+O(e^{-\eta})))=O(A{m}e^{A{m}}e^{-\eta}).$$
As a consequence, we may make the substitution
$${c}_{{m}}(s)=\sum_{{r}=1}^{{m}}\sum_{({p}_1,\dots,{p}_{{r}})\in P({m},{r})}\frac{(-1)^{{r}}}{{r}!}\prod_{{q}=1}^{{r}}\frac{1}{{p}_{{q}}}\sum_{w\in W^C_{{p}_{{q}}}}\left(e^{-sl(\tilde{\gamma}_w)}+O(A{p}_{{q}}e^{A{p}_{{q}}}e^{-\eta})\right)$$
\begin{equation}=\sum_{{r}=1}^{{m}}\sum_{({p}_1,\dots,{p}_{{r}})\in P({m},{r})}\frac{(-1)^{{r}}}{{r}!}\prod_{{q}=1}^{{r}}\frac{1}{{p}_{{q}}}\left(\left(\sum_{w\in W^C_{{p}_{{q}}}}e^{-sl(\tilde{\gamma}_w)}\right)+\left(\sum_{w\in W^C_{{p}_{{q}}}}O(A{p}_{{q}}e^{A{p}_{{q}}}e^{-\eta})\right)\right).\end{equation}
To handle this expression, first note that for any $\tilde{\gamma}_w$ of length ${m}$, {$l(\tilde{\gamma}_w)\leq CLm$, so by $|s|\leq\frac{A}{L}$,}
$$|e^{-sl(\tilde{\gamma}_w)}|\leq e^{{C}A{{m}}}.$$
Furthermore, $|W^C_{{p}_{{q}}}|\leq (2{d})^{{p}_{{q}}}\leq (2{d})^{2{d}}${, so} 
{\begin{equation}\label{maintermbound}\sum_{w\in W^C_{p_q}}e^{-sl(\tilde{\gamma}_w)}\leq (2d)^{2d}e^{CAm}\leq (2d)^{2d}e^{2CAd}=O_{g,n,C,A}(1).\end{equation}}
{By the same logic,
\begin{equation}\label{errortermbound}\sum_{w\in W^C_{p_q}}O(Ap_qe^{Ap_q}e^{-\eta})\leq O(A(2d)^{2d+1}e^{2Ad}e^{-\eta })=O(e^{-\eta })\end{equation}
Applying Equation \ref{errortermbound},
\begin{equation}c_m(s)=\sum_{r=1}^m\sum_{(p_1,...,p_r)\in P(m,r)}\frac{(-1)^r}{r!}\prod_{q=1}^r\frac{1}{p_q}\left(\left(\sum_{w\in W^C_{p_q}}e^{-sl(\tilde{\gamma}_w)}\right)+O(e^{-\eta })\right)\end{equation}
From here, expanding and applying Equation \ref{maintermbound},
$$\prod_{q=1}^r\frac{1}{p_q}\left(\left(\sum_{w\in W^C_{p_q}}e^{-sl(\tilde{\gamma}_w)}\right)+O(e^{-\eta })\right)$$
\begin{equation}\leq\left(\prod_{q=1}^r\sum_{w\in W^C_{p_q}}e^{-sl(\tilde{\gamma}_w)}\right)+\sum_{q'=1}^r\binom{r}{q'}O(1)^{r-q'}O(e^{-\eta })^{q'}\end{equation}
where in the second term the $\frac{1}{p_q}$ terms have been dropped, as this term will be an error term and all $\frac{1}{p_q}$ are at most $1$.
As $q'\geq 1$ in the latter sum and $r$ and $q'$ are bounded in terms of $g$ and $n$, this expression becomes
\begin{equation}\left(\prod_{q=1}^r\sum_{w\in W^C_{p_q}}e^{-sl(\tilde{\gamma}_w)}\right)+O(e^{-\eta }).\end{equation}
Therefore
\begin{equation}c_m(s)= \sum_{r=1}^m\sum_{(p_1,...,p_r)\in P(m,r)}\frac{(-1)^r}{r!}\left(\left(\prod_{q=1}^r\sum_{w\in W^C_{p_q}}e^{-sl(\tilde{\gamma}_w)}\right)+O(e^{-\eta })\right).\end{equation}
Recalling the definition of $\tilde{c}_m(s)$, this expression simplifies to
\begin{equation}c_m(s)=\tilde{c}_m(s)+\sum_{r=1}^m\sum_{(p_1,...,p_r)\in P(m,r)}\frac{(-1)^r}{r!}O(e^{-\eta}).\end{equation}
As $P(m,r)\leq (2d)^{2d}$, it is now straightforward to verify that
\begin{equation}\sum_{r=1}^m\sum_{(p_1,...,p_r)\in P(m,r)}\frac{(-1)^r}{r!}O_{g,n,C,A}(e^{-\eta})=O_{g,n,C,A}(e^{-\eta}).\end{equation}
We conclude that
\begin{equation}c_m(s)-\tilde{c}_m(s)|=O_{g,n,C,A}(e^{-\eta}).\ \square\end{equation}

\section{Proof of Main Results}

\textbf{Proof of Theorem \ref{mainresult}} Assume $\Gamma$ has ${d}$ edges, and recall that ${d}=O_{g,n}(1)$. Take cycle expansions and apply Lemmas \ref{truncation} and \ref{earlyconvergence}. Then:
$$|d(s,z)-\tilde{d}(s,z)|$$
$$=\left|\sum_{{m}=1}^\infty z^{{m}} {c}_{{m}}(s)-\sum_{{m}=1}^{2{m}}z^{{m}}\tilde{{c}}_{{m}}(s)\right|$$
$$\leq\sum_{{m}=1}^{2{d}}|z|^{{m}}|{c}_{{m}}(s)-\tilde{{c}}_{{m}}(s)|+\left|\sum_{{m}=2{d}+1}^\infty z^{{m}}{c}_{{m}}(s)\right|$$
$$=2dA^{2d}O(e^{-{\eta}})+O(e^{-{(1-\epsilon)\eta}})$$
\begin{equation}=O_{g,n,C,A,{\epsilon}}(e^{-{(1-\epsilon)\eta}}).\text{ $\square$}\end{equation}

\

\textbf{Proof of Corollary \ref{maincorollary}.} Recall that $MRG_{g,n}$ has a natural cell complex structure, where cells correspond to topological ribbon graph isomorphism classes of ribbon graphs with genus $g$ and $n$ boundaries. The assumption that $\Gamma$ is trivalent further implies that $\Gamma$ lies in the interior of a top-dimensional cell, so for sufficiently large ${t}$, all $\frac{1}{\alpha_{{t}}}\Gamma_{{t}}$ have the same topological type as $\Gamma$.  Furthermore, {the topology on $MRG_{g,n}$ implies that} the edge lengths of these $\frac{1}{\alpha_{{t}}}\Gamma_{{t}}$ converge to the lengths of the corresponding edges of $\Gamma$.  {Formally, let $\epsilon_1,...,\epsilon_{2d}$ be the directed edges of $\Gamma$ (as in Subsection \ref{zetafunctions}), and for sufficiently large $t$, let $\epsilon_{1,t},...,\epsilon_{2d,t}$ be the corresponding edges of $\Gamma_t$. Then for each $\epsilon_j$, $\lim_{t\to\infty}\frac{l(\epsilon_{j,t})}{\alpha_t}\to l(\epsilon_j)$. As a consequence, if $\Gamma_t\in MRG_{g,n}(L_{1,t},...,L_{n,t})$, then each sequence $(L_{h,t})$ satisfies $\lim_{t\to\infty}\frac{L_{h,t}}{\alpha_t}\to L_h$.}

{Therefore, there is some $T_0$ so that for any fixed $X_t$ with $t\geq T_0$, the $L_{h,t}$ are, say, bounded above by $2\max(L_1,...,L_n)\alpha_t$. Furthermore, if $E$ is defined as the minimum edge length of $\Gamma$, then there exists a $T_1$ so that for any fixed $X_t$ with $t\geq T_1$, every edge of $\Gamma_t$ has length at least $\frac{E}{2}\alpha_t$, and it follows that $X_t$ is $\frac{E}{2}\alpha_t$-great.}

{Recall that $U$ is a bounded domain by assumption; we will restrict to points $(s,z)\in\complex^2$ with $s\in U$ and $|z|\leq 2$. Let $J_s$ be the matrix derived from $\Gamma$ as in Subsection \ref{zetafunctions}, and let $J_{s,t}$ be the corresponding matrix for $\frac{1}{\alpha_t}\Gamma_t$. Then by the edge length limits given above and by Equation \ref{graphmatrixdef}, we have that $\lim_{t\to\infty} J_{s,t}\to J_s$ uniformly for $s\in U$, where the matrix convergence is entrywise. From Equation \ref{graphzetadef}, this limit is sufficient to conclude that the dynamical zeta functions $\tilde{d}_{\frac{1}{\alpha_{{t}}}\Gamma_{{t}}}(s,z)$ converge uniformly in $s$ and $z$ to $\tilde{d}_{\Gamma}(s,z)$.}

{As $\alpha_{{t}}\to\infty$, Theorem \ref{mainresult} and the fact that $X_t$ is $\frac{E}{2}\alpha_t$-great} implies that $\sup_{|z|\leq 2,s\in \frac{1}{\alpha_{{t}}}U}|d_{X_{{t}}}(z,s)-\tilde{d}_{\Gamma_{{t}}}(z,s)|$ converges to $0$ as ${t}\to\infty$. Since $\tilde{d}_{\Gamma_{{t}}}(z,s)$ and $\tilde{d}_{\frac{1}{\alpha_{{t}}}\Gamma_{{t}}}(z,s)$ are related by the easily verified formula
\begin{equation}\tilde{d}_{\Gamma_{{t}}}(z,\alpha_{{t}}s)=\tilde{d}_{\frac{1}{\alpha_{{t}} }\Gamma_{{t}}}(z,s)\end{equation}
we may conclude that
\begin{equation}\lim_{{t}\to\infty}\sup_{|z|\leq 2, s\in U}|d_{X_{{t}}}(z,\alpha_{{t}}s)-\tilde{d}_{\Gamma_{{t}}}(z,s)|=0\end{equation}
Recalling the definition of $\Res(X_{{t}},\alpha_{{t}})$, the desired result now follows immediately from an application of Rouch{\'e}'s Theorem. $\square$

\printbibliography

\end{document}